\RequirePackage{fix-cm}
\documentclass[smallextended]{svjour3}       
\smartqed  
\usepackage{amsmath}
\usepackage{booktabs}
\usepackage{mathtools}
\usepackage{amssymb,epsfig,multirow}
\usepackage{amsbsy}
\usepackage{appendix}
\usepackage{graphicx}
\usepackage{adjustbox}
\usepackage{url}
\usepackage{float}
\usepackage{color}
\usepackage{enumerate}
\usepackage[short]{optidef}
\usepackage{algorithm}
\usepackage{algpseudocode}
\usepackage{comment}
\usepackage[symbol]{footmisc}

\begin{document}

\title{Fix and Bound: An efficient approach for solving large-scale quadratic programming problems with box constraints}

\titlerunning{Fix and Bound: An efficient approach for solving large-scale BoxQPs}       

\author{Marco Locatelli \and Veronica Piccialli\footnote{Corresponding author.} \and \\ Antonio M. Sudoso}


\institute{Marco Locatelli \at
              Dipartimento di Ingegneria e Architettura, Universit\`{a} di
Parma, Parco Area delle Scienze, 181/A, Parma, 43124, Italy \\
              \email{marco.locatelli@unipr.it}
           \and
           Veronica Piccialli \at
              Dipartimento di Ingegneria Informatica, Automatica e Gestionale ``A. Ruberti'', Sapienza Universit\`{a} di Roma, Via Ariosto 25, Rome, 00185, Italy \\
              \email{veronica.piccialli@uniroma1.it}
           \and
           Antonio M. Sudoso \at
            Dipartimento di Ingegneria Informatica, Automatica e Gestionale ``A. Ruberti'', Sapienza Universit\`{a} di Roma, Via Ariosto 25, Rome, 00185, Italy \\
              \email{antoniomaria.sudoso@uniroma1.it}
}

\date{Received: date / Accepted: date}

\maketitle

\begin{abstract}
In this paper, we propose a branch-and-bound algorithm for solving nonconvex quadratic programming problems with box constraints (BoxQP). Our approach combines existing tools, such as semidefinite programming (SDP) bounds strengthened through valid inequalities, with a new class of optimality-based linear cuts which leads to variable fixing. The most important effect of fixing the value of some variables is the size reduction along the branch-and-bound tree, allowing to compute bounds by solving SDPs of smaller dimension. Extensive computational experiments over large dimensional (up to $n=200$) test instances show that our method is the state-of-the-art solver on large-scale BoxQPs. Furthermore, we test the proposed approach on the class of binary QP problems, where it exhibits competitive performance with state-of-the-art solvers.
\end{abstract}

\section{Introduction}
\label{sec:intro}
In a nonconvex QP problem a quadratic function is minimized over a polytope {$\mathcal{P}$}, i.e.:
\begin{equation*}
\begin{aligned}
\min_{} \quad & \frac{1}{2} {\bf x }^\top {\bf Q} {\bf x} + {\bf c}^\top {\bf x}  \\
\textrm{s.t.} \quad & {\bf x}\in {{\cal P}},\\
\end{aligned}
\end{equation*}
where ${\bf Q}\in \mathbb{R}^{n\times n}$ is symmetric and not positive semidefinite, ${\bf c}\in \mathbb{R}^n$.
The class of QP problems include as special cases Standard QP (StQP) problems, where the polytope {$\mathcal{P}$} is the unit simplex, and BoxQP problems, where the polytope {$\mathcal{P}$} is the unit box. 
Both special cases are NP-hard. Indeed, in \cite{Motzkin65} a polynomial reduction of Max-Clique problems into StQP problems is provided, while Max-CUT problems can be reformulated as BoxQP problems. 
Despite their difficulties, QP problems have a significant structure that allows the definition of different branch-and-bound (B\&B) approaches. For this reason, they attracted considerable attention in the literature (see, e.g., \cite{Bonami18,Burer09a,Chen12,Gondzio18,Liuzzi19,Liuzzi22,Xia20}), and even well-known commercial solvers such as {\tt Cplex} and {\tt Gurobi} offer exact approaches for their solution. 

QP problems can be reformulated in different ways. In the bilinear reformulation, each linear function ${\bf Q}_i {\bf x}$ is set equal to a new variable $y_i$ and substituted by it in the objective function (here ${\bf Q}_i$ denotes the $i$-th row of matrix ${\bf Q}$). In spite of its simplicity, the bilinear reformulation turned out to be the best one for a special class of QP problems arising from game theory \cite{Liuzzi20}. 
In the spectral reformulation (see, e.g., \cite{Sahinidis21}), first a spectral decomposition of matrix ${\bf Q}$ is performed, then each term ${\bf u}_i^\top {\bf x}$, where ${\bf u}_i$ is an eigenvector associated to a negative eigenvalue of ${\bf Q}$,  is set equal to a variable $z_i$ and replaced into the objective by such variable. In \cite{Liuzzi22} the spectral reformulation is employed when the number of negative eigenvalues of matrix ${\bf Q}$ is below a given threshold, while, alternatively, the bilinear reformulation is employed. This approach turns out to be quite robust and highly competitive with existing solvers (including {\tt Gurobi} and  {\tt Cplex}) over general QP problems. In \cite{Burer09} a copositive reformulation of QP problems has been introduced, where the problem is reformulated as a linear problem over the completely positive cone. The result in \cite{Burer09} extends previous results, in particular the result in  \cite{Bomze02}, where it is shown that StQP problems can be reformulated as linear problems over the copositive cone (the dual of the completely positive cone). The copositive reformulation is interesting from the theoretical point of view since it shows that nonconvex QP problems can be reformulated as convex problems. Unfortunately, the cone of completely positive matrices, though convex, is not tractable, so that existence of a convex reformulation does not lead to a polynomial algorithm. However, the copositive reformulation can be exploited to define convex relaxations, by replacing the not tractable cone of completely positive matrices with tractable cones such as the semidefinite or nonnegative one. 
In the KKT reformulation (see, e.g., \cite{Hansen93}), additional variables and constraints are introduced, the former corresponding to Lagrange multipliers of the constraints, the latter to the KKT conditions, and, in addition, the objective function is linearized by exploiting the KKT conditions. Strictly related to the KKT reformulation is the Mixed Integer Linear Programming (MILP) reformulation (see, e.g., \cite{Xia20}), where the complementarity conditions are linearized through the introduction of binary variables.

All these reformulations allow to define different convex relaxations (linear, convex quadratic, semidefinite relaxations) with different computational costs but also with a different quality of the corresponding lower bounds, thus posing the question about which is the best compromise between quality and computational cost. In fact, it is observed in practice that the best choice is strictly related to the problem at hand, and some works (see, e.g., \cite{Liuzzi22,Sahinidis21}) suggest to adaptively choose between different relaxations depending on specific properties of the problem to be solved, possibly making different choices at different nodes of the B\&B tree. 

Relaxations can be strengthened through the addition of valid inequalities. These include, e.g., Reformulation Linearization Technique (RLT) inequalities and triangle inequalities which will be discussed in what follows. A special class of valid inequalities are those involving variable bounds. Procedures to define such inequalities are known as bound-tightening techniques  (see, e.g., \cite{Caprara10,Gleixner17,Sahinidis96,Tawarmalani04}). In the context of QP problems, optimality-based bound tightening techniques are employed, i.e., such techniques define new bounds for the variables which may not be fulfilled by some feasible solutions but are always fulfilled by optimal solutions. Although they might be quite expensive, they often turn out to be of primary importance to enhance the performance of B\&B approaches. We refer to \cite{Gleixner17} for some strategies to find a good compromise between the effectiveness and cost of bound-tightening techniques.
{\section*{Statement of contribution}
In this paper, we propose a B\&B algorithm for solving Box QP problems. The approach combines tools already known in the literature, such as SDP-based bounds strengthened by adding valid inequalities (in particular, RLT and triangle inequalities) with a newly proposed variable fixing procedure. This procedure has two features, which represent the main original contributions of our work:
\begin{itemize}
\item The procedure is a {\em multiple} variable fixing procedure, i.e., the values of multiple variables are fixed at the same time; while {\em single} variable fixing procedures have been proposed in the literature (see, e.g., \cite{galli2018binarisation}, where variable fixing is employed within the framework of a heuristic approach for Box QP problems), we are not aware of multiple variable fixing procedures;
 \item The procedure is based on the solution of SDP problems; in the literature SDPs are usually solved to compute bounds, while we are not aware of SDP-based variable fixing procedures.
 \end{itemize}
Solving SDPs is expensive and the fixing procedure adds new SDPs to be solved. However, the strong SDP relaxation obtained by adding valid inequalities makes the multiple fixing procedure effective, allowing for a reduction of the size of the subproblems solved at nodes of the B\&B tree. Once variables are fixed at a node, they are also fixed in all its successors, thus reducing the time both for the new bound computations and for the further variable fixing procedures. The experiments show
 that, over the most challenging instances of Box QP and Binary QP problems, the inclusion of this variable fixing procedure leads to a reduction of the computing times ranging between 20\% and 50\%.
 Further experiments also show that the proposed
 solution approach is highly competitive with the existing solvers.}

The remainder of the paper is structured as follows.
In Section \ref{sec:boxqp} we introduce the class of BoxQP problems and discuss different features which can be exploited for solving them more efficiently.
In Section \ref{sec:sdpbound} we will discuss a bound based on the solution of a semidefinite problem obtained by strengthening the classical Shor relaxation with RLT and triangle inequalities. Moreover, we discuss the use of first-order SDP solvers, which allows to compute SDP bounds more efficiently with respect to the previous generation of SDP solvers.
In Section  \ref{sec:fixing} we discuss a class of optimality-based linear cuts which, in particular, leads to variable fixing, i.e., the possibility of fixing the value of some variables, thus allowing for a dimension reduction.
In Section \ref{sec:standard} we briefly present further components which are needed in order to have a full implementation of a B\&B approach. Such components are defined in a standard way, taken from the literature. 
Finally, in Section \ref{sec:exp} we present and discuss different computational experiments, whereas in Section \ref{sec:conclusions} we draw some conclusions and highlight possible future research directions.

\subsection{Notation}
Throughout this paper we use the following notation: $\mathcal{S}^n$ denotes the space of $n\times n$ real symmetric matrices, $\mathbb{R}^n$ is the space of $n$-dimensional real vectors and $\mathbb{R}^{m \times n}$ is the space of $m \times n$ real matrices.
We denote by $\mathbf{X} \succeq {\bf O}$ a matrix $\mathbf{X}$ that is positive semidefinite and by $\mathcal{S}^n_+$ be the set of positive semidefinite matrices of size $n \times n$. We denote by $\mathbf{0}^n, \mathbf{1}^n$ the $n$-dimensional vectors of all zeros and of all ones. We omit the superscript in case the dimension is clear from the context. Given $\mathbf{x} \in \mathbb{R}^n$, $\textrm{Diag}(\mathbf{x})$ is the $n \times n$ diagonal matrix with $\mathbf{x}$ on its
diagonal; and given $\mathbf{X} \in \mathbb{R}^{n \times n}$, $\textrm{diag}(\mathbf{X})$ is the vector with the diagonal elements of $\mathbf{X}$. We denote by $\bullet$ the 
trace inner product. That is, for any
$\mathbf{A}, \mathbf{B} \in \mathbb{R}^{m\times n}$, we define $\mathbf{A} \bullet \mathbf{B} := \textrm{trace}(\mathbf{B}^\top \mathbf{A})$. For a vector $\mathbf{v}$ and an index set $I$, $\mathbf{v}^I$ is defined as the
vector composed of entries of $\mathbf{v}$ that are indexed by $I$. Also, given a matrix $\mathbf{X} \in \mathbb{R}^{n \times n}$, $\mathbf{X}^{II}$ is defined as the matrix composed of entries of $\mathbf{X}$ whose rows and columns
are indexed by $I$. 

\section{BoxQP problems}
\label{sec:boxqp}
In BoxQP problems the feasible polytope {$\mathcal{P}$} is the $n$-dimensional unit box $\{{\bf x}\ :\ {\bf 0}^n\leq {\bf x}\leq {\bf 1}^n\}$. Note that the feasible region of a BoxQP problem can also be a generic box but, obviously, any box can be replaced by the unit box, possibly after scaling and translating the variables.
BoxQP problems have some specific properties which can be exploited. 
Let $N=\{i\ :\ Q_{ii}\leq 0\}$ and $P=\{i\ :\ Q_{ii}>0\}$.
First, we notice (see {Proposition 1 in } \cite{Hansen93}) that: 
$$
i\in N \Rightarrow x_i \in \{0,1\},
$$
at optimal solutions of Box QP problems, so that some of the variables can be immediately considered as binary ones. 
Next, we observe that the stationarity condition in the KKT system reduces to:
$$
{\bf Q}{\bf x}+{\bf c}-\boldsymbol{\mu}+\boldsymbol{\gamma}={\bf 0},
$$
where $\boldsymbol{\mu}$ is the Lagrange multiplier vector of the nonnegativity constraints, while $\boldsymbol{\gamma}$ is the Lagrange multiplier vector of constraints ${\bf x}\leq \mathbf{1}^n$.
Then, for all $i\in P$, we have that:
\begin{equation}
\label{eq:kktboxqp}
\begin{array}{lllll}
x_i=0 & \Rightarrow & \gamma_i=0 & \Rightarrow & {\bf Q}_i {\bf x}\geq - c_i \\ [6pt]
x_i=1 & \Rightarrow & \mu_i=0 & \Rightarrow & {\bf Q}_i {\bf x}\leq - c_i \\ [6pt]
0<x_i<1 & \Rightarrow & \gamma_i=0,\ \mu_i=0 & \Rightarrow & {\bf Q}_i {\bf x}= - c_i,
\end{array}
\end{equation}
while for all $i\in N$, one of the first two implications holds.
In B\&B approaches for general QPs, spatial branching is often adopted, i.e., first a variable is selected, together with a threshold value for it, then the value of such variable is imposed to be not larger than the threshold in one child node, and not smaller than the threshold in the other child node. But, in view of the above properties, in BoxQPs spatial branching can be replaced by binary branching, if the selected variable $x_i$ is such that $i\in N$, or ternary branching, if the selected variable $x_i$ is such that $i\in P$. In the latter case in each one of the three child nodes one of the conditions $\{x_i=0\}$, $\{x_i=1\}$, and $\{{\bf Q}_i {\bf x}=-c_i\}$ is imposed.
{While spatial branching only allows to prove convergence to 0 of the difference between the upper and lower bound values computed by the B\&B approach, binary/ternary branching has the theoretical advantage that such difference becomes equal to 0 after the generation of a finite, though possibly exponential, number of B\&B nodes.}\newline
A further way to exploit the properties of variables $x_i$ for $i\in N$ and $i\in P$, which also represents the main contribution of this work, will be thoroughly discussed in Section \ref{sec:fixing}.
The same properties also allow to introduce a MILP reformulation of BoxQPs. All variables $x_i$, $i\in N$, are immediately set as binary.
For each $i\in P$, we first denote by $\ell_i$ and $u_i$ a lower and an upper bound, respectively, for ${\bf Q}_i {\bf x}$ over the unit box.
In particular, we can simply set:
$$
\ell_i=\sum_{j\ :\ Q_{ij}<0} Q_{ij},\ \ \ u_i=\sum_{j\ :\ Q_{ij}>0} Q_{ij}.
$$
Next, again for each $i\in P$ we introduce three binary variables $\delta_i, \rho_i, \xi_i$, corresponding to $x_i=0$, $x_i=1$ and $0<x_i<1$, respectively.
Thus, in order to formulate (\ref{eq:kktboxqp}), we need the following constraints for each $i\in P$:
$$
\begin{array}{l}
x_i\leq (1-\delta_i),\ \ \ {\bf Q}_i{\bf x}\geq \ell_i(1-\delta_i) -\delta_i c_i  \\ [6pt]
x_i\geq \rho_i,\ \ \ {\bf Q}_i{\bf x}\leq u_i(1-\rho_i) -\rho_i c_i  \\ [6pt]
{\bf Q}_i{\bf x}\leq u_i(1-\xi_i) -\xi_i c_i,\ \ \ {\bf Q}_i{\bf x}\geq \ell_i(1-\xi_i) -\xi_i c_i \\ [6pt]
\delta_i+\rho_i+\xi_i=1.
\end{array}
$$
However, according to our personal experience {and} according to the
experiments reported in \cite{Xia20}, solution approaches based on MILP reformulations are not quite efficient for BoxQPs, while they are for another special class of QP problems, namely the already mentioned class of StQP problems (see \cite{Gondzio18,Xia20}). 
\newline\newline\noindent
It is also worthwhile to remark that any approach for the solution of BoxQPs is also a valid approach for the solution of binary QPs, i.e.,
problems
\begin{equation}\label{prob:binaryQP}
  \begin{array}{lll}
\min_{{\bf x}\in \{0,1\}^n} & \frac{1}{2} {\bf x }^\top {\bf Q} {\bf x} + {\bf c}^\top {\bf x}.& 
\end{array}  
\end{equation}
Indeed, taking into account that the binary condition for a variable $x_i$ can also be written as $x_i(1-x_i)=0$, we can convert the binary QP problem into the BoxQP problem
\begin{equation}\label{prob:BoxQPLambda}
    \begin{array}{lll}
\min_{{\bf x}\in [0,1]^n} & \frac{1}{2} {\bf x }^\top ({\bf Q}-\boldsymbol{\Lambda}) {\bf x} + ({\bf c}+{\frac{1}{2}}\boldsymbol{\lambda})^\top {\bf x},& 
\end{array}
\end{equation}
where $\boldsymbol{\Lambda}=\textrm{Diag}(\boldsymbol{\lambda})$ and for each $i\in \{1,\ldots,n\}$ we have $-\lambda_i+Q_{ii}\leq 0$, i.e., the Hessian matrix of the BoxQP problem has nonpositive diagonal elements, so that, by the previously mentioned
property, optimal solutions of the BoxQP problem are forced to be binary solutions and, thus, to be optimal solutions also of the Binary QP problem.
Possible choices for the vector $\boldsymbol{\lambda}$ will be discussed later on.

In the following, we introduce a B\&B approach for the solution of BoxQPs. First, we present an SDP-based technique to compute lower bounds at nodes of the B\&B tree, with the addition of valid inequalities taken from the literature. Moreover, we discuss the use of first-order SDP solvers for a more efficient computation of the bounds. Next, we present a novel technique which allows to fix at the same time different variables to 0 or to 1 through the solution of an auxiliary SDP problem. Such novel technique is proven to be quite beneficial at large dimensions, as we show in the computational experiments.
\section{SDP bounds for BoxQP with valid inequalities}
\label{sec:sdpbound}
Relaxation techniques for BoxQP are mainly based on linearization, convex quadratic programming, or semidefinite programming. To compute such a relaxation, the quadratic objective function is {typically} reformulated as a linear function in an extended space of variables. More precisely, a symmetric matrix variable $\mathbf{X}$ satisfying $\mathbf{X}=\mathbf{x}\mathbf{x}^\top$ is introduced
and the equivalent formulation in the $(\bf{x}, \bf{X})$-space is then solved by means of a B\&B algorithm based on a
relaxation of the latter non-convex equalities. The most commonly used relaxation is the so-called McCormick relaxation which is obtained by replacing the nonlinear equalities with RLT constraints \cite{mccormick1976computability}. On the other hand, the semidefinite relaxation due to Shor \cite{shor1987quadratic} can be obtained by relaxing $\mathbf{X}=\mathbf{x}\mathbf{x}^\top$ to $\mathbf{X} - \mathbf{x}\mathbf{x}^\top \succeq 0$ which in turn is equivalent to the matrix inequality
\begin{equation*}
    \mathbf{\tilde{X}} = \begin{bmatrix} 
	1 & \mathbf{x}^\top\\
	\mathbf{x} & \mathbf{X}\\
	\end{bmatrix} \succeq {\bf O},
\end{equation*}
due to the Schur complement. Combining the classical Shor relaxation with the RLT inequalities leads to strong bounds, as shown in \cite{anstreicher2009semidefinite}. Therefore, the basic SDP relaxation we use in our bounding routine is
\begin{equation}
\begin{aligned}\label{eq:SDPbound}
\min_{} \quad & \frac{1}{2} \bf{Q} \bullet \bf{X} + {\bf c}^\top {\bf x}  \\
\textrm{s.t.} \quad & {\bf 0}\le {\bf x}\le \mathbf{1}\\
& \bf{X}-{\bf x}{\bf x}^\top \succeq  {\bf O} \\
& (\bf{x}, \bf{X}) \in \mathcal{S}_{\textrm{RLT}}\\
\end{aligned}
\end{equation}
where, for $1 \leq i < j \leq n$, the set $\mathcal{S}_{\textrm{RLT}} \subset \mathbb{R}^n \times \mathbb{S}^n$ is defined as
\begin{align*}
    \mathcal{S_{\textrm{RLT}}} = \left\{(\mathbf{x}, \mathbf{X}) \in \mathbb{R}^n \times \mathbb{S}^n : \begin{array}{l}
    X_{ij} \geq 0, \ X_{ij} \leq x_i, \ X_{ij} \leq x_j, \ X_{ij} \geq x_i + x_j - 1 \\
    \end{array} \right\}.
\end{align*}
This kind of SDP is a doubly nonnegative program (DNN) since the matrix variable $\mathbf{\tilde{X}}$ is both positive semidefinite and element-wise nonnegative. Note that the original bound constraints ${\bf 0}^n \le {\bf x}\le \mathbf{1}^n$ are implied by the RLT constraints and could be removed (see \cite{sherali1995reformulation}, Proposition 1).
In order to derive stronger relaxations in the $(\bf{x}, \bf{X})$-space, Yajima \& Fujie \cite{Yajima1998} studied the convex hull of feasible solutions of the BoxQP problem, i.e., the set
\begin{align*}
        \textrm{QPB}_n = \textrm{conv}\left\{ (\mathbf{x}, \mathbf{X}) \in [0, 1]^{{n+}\binom{n+1}{2}} : \ X_{ij} = x_i x_j \ (1 \leq i \leq j \leq n)\right\},
\end{align*}
and they showed valid inequalities for it. 
In \cite{burer2009nonconvex}, Burer \& Letchford studied $\textrm{QPB}_n$ in more detail and they showed that any inequality valid for the boolean quadratic polytope of order $n$ (see \cite{padberg1989boolean}), i.e., the set 
\begin{align*}
        \textrm{BQP}_n = \textrm{conv}\left\{ (\mathbf{x}, \mathbf{X}) \in \{0, 1\}^{{n+}\binom{n}{2}} : \ X_{ij} = x_i x_j \ (1 \leq i < j \leq n)\right\},
\end{align*}
is valid also for $\textrm{QPB}_n$. One important class of these cuts are the triangle inequalities defined as follows:
\begin{equation}
\begin{split}
\label{eq:triangle}
x_i + x_j + x_t - X_{ij} - X_{it} - X_{jt} &\leq 1\\
X_{ij}+X_{it}-x_i-X_{jt} &\leq 0 \\
X_{ij}+X_{jt}-x_j-X_{it} &\leq 0 \\
X_{it}+X_{jt}-x_t-X_{ij} &\leq 0.
\end{split}
\end{equation}
These hold for all distinct triplets $(i, j, t)$. In order to obtain tight lower bounds, we use the cutting-plane approach where triangle inequalities are iteratively added and purged after each bound computation. The enumeration of all triangle inequalities is computationally inexpensive, even for large instances. However, adding all of them would make the SDP relaxation intractable. To keep the SDP to a modest size, we limit the number of triangle inequalities that we can add at each cutting-plane iteration. First, the optimal solution of the basic SDP relaxation in \eqref{eq:SDPbound} is computed. Then, we separate the initial set of triangle inequalities, we sort them by the magnitude of violation and we add the most violated ones {up to a maximum number}. Next, we remove all inactive constraints and we add new violated triangles. Finally, the problem with an updated set of inequalities is solved and the process is repeated as long as the increase of the lower bound is sufficiently large. Indeed, we terminate the bounding routine when the relative difference between consecutive lower bounds is less than a fixed threshold. If the gap cannot be closed after some cutting-plane iterations, we terminate the bounding routine and we branch. Note that, to improve the efficiency of the overall B\&B algorithm, we inherit {triangle inequalities} from parent to children nodes. This allows us to save a significant number of cutting-plane iterations, and therefore computational time, when processing subproblems.
\subsection{First-order SDP solver with reliable bounds}
\label{sec:firstorder}
The high quality of SDP bounds with the addition of triangle inequalities allows to reach small gaps between upper and lower bound even at the root of the B\&B tree (see also the experiments in Section \ref{sec:moseksdpnal}). The limitation is that SDP solvers do not scale well as the problem size increases. 
In particular, interior-point methods (IPMs) can solve SDPs to arbitrary accuracy in polynomial time, but they are not capable of solving large SDPs with a large number of inequalities, including nonnegativity constraints \cite{alizadeh1995interior}.   
Compared to IPMs, first-order solvers based on augmented Lagrangian methods (ALMs) can scale to significantly larger problem sizes, while trading off the accuracy of the resulting output \cite{sun2015convergent,yang2015sdpnal}. 
However, these methods require some extra care when embedded into a B\&B algorithm. 
When using first-order methods, it is hard to reach a solution to high precision in a reasonable amount of time. For this reason, we consider post-processing methods to guarantee a valid lower bound for the SDP relaxation. This lower bound is a safe underestimate for the optimal solution of the SDP relaxation that can in turn be used as a lower bound on the original BoxQP problem. 

Duality results in semidefinite programming imply that the objective function value of every feasible solution of the dual problem yields a bound on the optimal objective function value of the primal. Every dual feasible solution and, in particular, the optimal dual solution of an SDP relaxation, gives a valid bound on the solution of the original optimization problem. However, the dual objective function value represents a valid dual bound only if the SDP relaxation is solved to high precision. Two different post-processing techniques can be defined: one that adds a negative perturbation to the dual objective function value \cite{jansson2008rigorous} and one that generates a dual feasible solution, and hence a bound, by solving a linear program \cite{cerulli2021improving}. The first method is effective when a suitable bound on the maximum eigenvalue of an optimal solution of the primal SDP is known, as for example in \cite{piccialli2022a,piccialli2022b,piccialli2022c}. When dealing with BoxQP problems, instead, no tight bound can be derived, so the LP-based post-processing is a better option. 
First of all, we derive the dual problem of the SDP relaxation in \eqref{eq:SDPbound}.
Let $\tilde{\mathbf{Q}} = \frac{1}{2}\begin{bmatrix}
0 & \mathbf{c}^\top\\
\mathbf{c} & \mathbf{Q}
\end{bmatrix}$ and $\tilde{\mathbf{X}} = \begin{bmatrix}
x_0 & \mathbf{x}^\top\\
\mathbf{x} & \mathbf{X}
\end{bmatrix}$. 
We re-write problem \eqref{eq:SDPbound} as
\begin{equation}
\begin{aligned}\label{prob:SDPstandard}
\min_{} \quad & \tilde{\mathbf{Q}} \bullet \tilde{\mathbf{X}}  \\
\textrm{s.t.} \quad & x_0 = 1\\
& \mathcal{B}(\tilde{\mathbf{X}}) \geq \mathbf{b}\\
& \tilde{\mathbf{X}} \geq \mathbf{0}, \ \tilde{\mathbf{X}} \succeq  {\bf O},\\
\end{aligned}
\end{equation}
where $\mathcal{B} : \mathcal{S}^{n+1} \rightarrow \mathbb{R}^{m}$ defined by $\mathcal{B}(\tilde{\mathbf{X}}):=[\tilde{\mathbf{X}} \bullet \mathbf{B}^1, \dots, \tilde{\mathbf{X}} \bullet \mathbf{B}^m]^\top$ with $\mathbf{B}^i \in \mathcal{S}^{n+1}$, is the linear operator representing RLT and triangle inequalities on $\mathbf{\tilde{X}}$, and $\mathbf{b} \in \mathbb{R}^m$ is the right-hand side vector.
Consider multipliers $y \in \mathbb{R}$, $\mathbf{v} \in \mathbb{R}^m$, $\mathbf{S}$, $\mathbf{Z} \in \mathcal{S}^{n+1}$. The dual of problem \eqref{prob:SDPstandard} can be written as
\begin{equation}
\begin{aligned}\label{prob:SDPstandardDual}
\max_{} \quad & y + \mathbf{b}^\top \mathbf{v}  \\
\textrm{s.t.} \quad & \mathbf{E}_{11} y + \mathcal{B}^\star(\mathbf{v}) + \mathbf{Z} + \mathbf{S} = \tilde{\mathbf{Q}}\\
& \mathbf{S} \geq \mathbf{0}, \ \mathbf{v} \geq \mathbf{0}, \ \mathbf{Z} \succeq  {\bf O},\\
\end{aligned}
\end{equation}
where $\mathcal{B}^\star : \mathbb{R}^{m}  \rightarrow \mathcal{S}^{n+1}$ is the adjoint operator to $\mathcal{B}$ defined by $\mathcal{B}^\star(\mathbf{v}) := \sum_{i=1}^{m} v_i \mathbf{B}^i$ and $\mathbf{E}_{11}$ is the matrix of all zeros except entry $(1,1)$, which is equal to 1.
Let $({y}, {\mathbf{v}}, {\mathbf{S}},  {\mathbf{Z}})$ be any solution to the dual SDP \eqref{prob:SDPstandardDual} obtained by running an ALM method up to a certain precision. Note that this solution may be close to optimal solution but not necessarily dual feasible.
From the given $\mathbf{Z}$ we obtain the new positive semidefinite matrix $\tilde{\mathbf{Z}}$
by projecting $\mathbf{Z}$ onto the positive semidefinite cone. 
{This process involves computing the eigenvalue decomposition of $\mathbf{Z}$ and setting all negative eigenvalues to zero, while keeping the non-negative eigenvalues unchanged}.
Then, we solve the linear program
\begin{equation}
\begin{aligned}\label{prob:LPstandardDual}
\max_{} \quad & y + \mathbf{b}^\top \mathbf{v}  \\
\textrm{s.t.} \quad & \mathbf{E}_{11} y + \mathcal{B}^\star(\mathbf{v}) + \mathbf{S} = \tilde{\mathbf{Q}} - \tilde{\mathbf{Z}}\\
& \mathbf{S} \geq \mathbf{0}, \ \mathbf{v} \geq \mathbf{0},\\
\end{aligned}
\end{equation}
with respect to $(y, \mathbf{v}, \mathbf{S})$. If problem \eqref{prob:LPstandardDual} admits an optimal solution $(\tilde{y}, \tilde{\mathbf{v}}, \tilde{\mathbf{S})}$, then $(\tilde{y}, \tilde{\mathbf{v}}, \tilde{\mathbf{S}}, \tilde{\mathbf{Z}})$ is a feasible solution for problem \eqref{prob:SDPstandardDual} and we can return a valid dual bound. If it is infeasible, we are neither able to construct a feasible dual solution nor to construct a dual bound. In this case, one should continue running the ALM to a higher precision and apply the procedure to the improved point. However, we point out that in our computational experiments we never needed to increase the precision of the SDP solver.

\section{A new class of optimality-based linear inequalities and variable fixing}
\label{sec:fixing}
Let us denote by ${\cal S}$ the feasible region of the SDP relaxation strengthened with RLT and triangle inequalities.
Moreover, let $UB$ be the current upper bound for the optimal value of BoxQP.
The region
$$
{\cal S}^\star={\cal S}\cap \left\{\frac{1}{2}{\bf Q}\bullet {\bf X}+{\bf c}^\top {\bf x}\leq UB\right\},
$$
may exclude some feasible solutions of BoxQP but contains the whole set of its optimal solutions (more precisely, the projection of such region over the space of variables ${\bf x}$ contains the set of optimal solutions).
Now, 
let $T\subseteq N$ and let $T_0, T_1$ be a partition of $T$. 
Let us solve the following SDP problem:
\begin{equation}
\label{eq:varfix}
\begin{array}{ll}
\ell_{T}=\min_{({\bf X},{\bf x})\in {\cal S}^\star} & \sum_{i\in T_0} (1-x_{i})+\sum_{i\in T_1} x_{i}.
\end{array}
\end{equation}
Recalling that for all variables in $N$ and, thus, in $T$, we can impose a binary condition,  the following is an optimality-based valid linear inequality:
\begin{equation}
\label{eq:optlincut}
 \sum_{i\in T_0} (1-x_{i})+\sum_{i\in T_1} x_{i}\geq \lceil \ell_T \rceil.
\end{equation}
In fact, in order to take into account numerical errors, the right-hand side of the inequality is replaced by $\lceil \ell_T -\varepsilon_1\rceil$ for some tolerance value {$\varepsilon_1\in (0,1)$, that we set in the experiments equal to $0.01$}.
The identification of such valid linear inequality requires the solution of an SDP problem and is, thus, quite expensive. But there is at least one situation when such high computational cost is worthwhile, namely when
$\lceil \ell_T -\varepsilon_1\rceil=\lvert T \rvert$. Indeed, in this case we are able to fix all variables in $T$, more precisely, we can fix to 0 all variables in $T_0$, and to 1 all variables in $T_1$.
As already mentioned in Section \ref{sec:intro}, solvers for QP problems highly benefit from the application of bound tightening procedures, in spite of their computational cost. The multiple fixing strategy, which is also quite expensive, can be viewed as a special case of bound tightening, where range of variables are reduced to singletons.
One advantage of multiple fixing is the improvement of the lower bounds, which may reduce the size of the B\&B trees. But, probably, the most important effect of multiple variable fixing is the size reduction. 
Let $F_0, F_1$ denote the set of indices of the variables fixed to 0 and 1, respectively, at the current node. Let $U=\{1,\ldots,n\}\setminus \{F_0 \cup F_1\}$ be the set of indices of all remaining variables.
Without loss of generality, we assume that the indices of the variables are ordered as follows $F_0, U, F_1$. We denote by $\mathbf{0}^{F_0}$ the zero vector of size $\lvert F_0\rvert$, $\mathbf{1}^{F_1}$ the unit vector of size
$\lvert F_1\rvert$, and by $\mathbf{x}^{U}$ the vector of not fixed variables.
Hence, the BoxQP objective function  can be written as
\begin{align*}
& \frac{1}{2}
\begin{bmatrix}
\mathbf{0}^{F_0}\\
\mathbf{x}^{U}\\
\mathbf{1}^{F_1}
\end{bmatrix}^\top
\begin{bmatrix}
\mathbf{Q}^{F_0 F_0} & \mathbf{Q}^{F_0 U} & \mathbf{Q}^{F_0 F_1}\\
\mathbf{Q}^{U F_0} & \mathbf{Q}^{UU} & \mathbf{Q}^{U F_1}\\
\mathbf{Q}^{F_1 F_0} & \mathbf{Q}^{F_1 U} & \mathbf{Q}^{F_1 F_1}\\
\end{bmatrix}
\begin{bmatrix}
\mathbf{0}^{F_0}\\
\mathbf{x}^{U}\\
\mathbf{1}^{F_1}
\end{bmatrix}+
\begin{bmatrix}
\mathbf{c}^{F_0}\\
\mathbf{c}^{U}\\
\mathbf{c}^{F_1}
\end{bmatrix}^\top
\begin{bmatrix}
\mathbf{0}^{F_0}\\
\mathbf{x}^{U}\\
\mathbf{1}^{F_1}
\end{bmatrix} \\ 
& =
\frac{1}{2}
\begin{bmatrix}
\mathbf{x}^{U}\\
\mathbf{1}^{F_1}\\
\end{bmatrix}^\top
\begin{bmatrix}
\mathbf{Q}^{UU} & \mathbf{Q}^{UF_1}\\
(\mathbf{Q}^{UF_1})^\top & \mathbf{Q}^{F_1 F_1}\\
\end{bmatrix}
\begin{bmatrix}
\mathbf{x}^{U}\\
\mathbf{1}^{F_1}\\
\end{bmatrix} +
\begin{bmatrix}
\mathbf{c}^{U}\\
\mathbf{c}^{F_1}\\
\end{bmatrix}^\top
\begin{bmatrix}
\mathbf{x}^{U}\\
\mathbf{1}^{F_1}\\
\end{bmatrix}\\ 
& = \frac{1}{2} (\mathbf{x}^U)^\top \mathbf{Q}^{UU} \mathbf{x}^U + (\mathbf{x}^U)^\top \mathbf{Q}^{UF_1} \mathbf{1}^{F_1} + (\mathbf{c}^U)^\top \mathbf{x}^U +
\xi^{F_1,U}, 
\end{align*}
where $\xi^{F_1,U}= \frac{1}{2} (\mathbf{1}^{F_1})^\top \mathbf{Q}^{F_1 F_1} \mathbf{1}^{F_1}  + (\mathbf{c}^{{F_1}})^\top \mathbf{1}^{F_1}$.
Therefore, after removing the constant term 
$\xi^{F_1,U}$, 
the bounding problem becomes
\begin{equation}
\begin{aligned}\label{eq:SDPboundFixing}
\min_{} \quad & \frac{1}{2} \mathbf{\bar{Q}} \bullet \mathbf{X}^U + \mathbf{\bar{c}}^\top \mathbf{x}^U  \\
\textrm{s.t.} \quad & \mathbf{0}^{\lvert U\rvert} \le \mathbf{x}^U \le \mathbf{1}^{\lvert U\rvert}\\
& \mathbf{X}^U -{\bf x}^U({\bf x}^U)^\top \succeq  {\bf O}\\
& (\mathbf{x}^U, \mathbf{X}^U) \in \mathcal{S}_{\textrm{RLT}},
\end{aligned}
\end{equation}
where $\bar{\mathbf{Q}} = \mathbf{Q}^{UU}$ and $\bar{\mathbf{c}} = \mathbf{Q}^{U F_1} \mathbf{1}^{F_1} + \mathbf{c}^{U}$. 
The size reduction of the BoxQP problem from dimension $n$ to $\lvert U \rvert$, implies also the size reduction of its SDP relaxation and, since computing times of
SDP solvers tend to increase quickly with the size, such reduction allows to compute lower bounds more quickly.

An important question related to multiple variable fixing is how to choose $T$ and how to partition it into $T_0$ and $T_1$. The most natural strategy is to exploit
the solution of the SDP relaxation over the region ${\cal S}$. If we denote by $({\bf X}^\star,{\bf x}^\star)$ the optimal solution of such relaxation, we can set:
$$
T_0=\{i\in N\ :\ x_{i}^\star \leq \varepsilon_2\},\ \ \ T_1=\{i\in N\ :\ x_{i}^\star \geq 1-\varepsilon_2\},
$$
where $\varepsilon_2>0$: then, we put in $T_0$ all variables close to 0 in the optimal solution of the relaxation, and in $T_1$ all variables close to 1 in the optimal solution. 
\newline\newline\noindent
Finally, we point out that, rather than {going }through the solution of  (\ref{eq:varfix}), an alternative way to perform multiple variable fixing is through the solution of the following SDP problem:
\begin{equation}
\label{eq:varfix1}
\begin{aligned}
\alpha_T=\min_{({\bf X},{\bf x})\in {\cal S}} \quad & \frac{1}{2}{\bf Q}\bullet {\bf X}+{\bf c}^\top {\bf x}  \\
\textrm{s.t.} \quad & \sum_{i\in T_0} (1-x_{i})+\sum_{i\in T_1} x_{i}\leq \lvert T \rvert-1+\varepsilon_1.
\end{aligned}
\end{equation}
Indeed, if $\alpha_T\geq UB$, i.e., if the lower bound of the BoxQP problem over the intersection of ${\cal S}$ with the region defined by the linear inequality $\sum_{i\in T_0} (1-x_{i})+\sum_{i\in T_1} x_{i}\leq \lvert T \rvert-1+\varepsilon_1$, is not smaller than the current upper bound, then {no solution with function value lower than $UB$ can be found in this region, so that we can restrict the attention to the region defined by the inequality $\sum_{i\in T_0} (1-x_{i})+\sum_{i\in T_1} x_{i}> \lvert T \rvert-1+\varepsilon_1$. But in such case, due to the binary conditions, we can fix all variables in $T$.} In fact, the solution of (\ref{eq:varfix1}) is slightly less informative with respect to the solution of (\ref{eq:varfix}). Indeed, in case fixing is not possible, i.e., when
$\lceil \ell_T -\varepsilon_1\rceil<\lvert T \rvert$ concerning (\ref{eq:varfix}), or $\alpha_T<UB$ concerning (\ref{eq:varfix1}), the solution of (\ref{eq:varfix}) allows at least to add the linear cut (\ref{eq:optlincut}).
On the other hand, recent first-order SDP solvers discussed in Section \ref{sec:firstorder} encounter some difficulties when solving
(\ref{eq:varfix}), while they are able to solve (\ref{eq:varfix1}) quite efficiently. For this reason, in our experiments we adopted (\ref{eq:varfix1}). {Note that Problem (\ref{eq:varfix1}) always admits feasible solutions. Indeed, any solution of the form $(\mathbf{x}, \mathbf{X})$ where $\mathbf{x}$ is a binary vector satisfying the additional inequality in (\ref{eq:varfix1}), and $\mathbf{X}=\mathbf{x}\mathbf{x}^\top$, is feasible.}
\noindent
We summarize the fixing procedure in Algorithm \ref{algo:fixing}. The number of variables that can be fixed by Algorithm \ref{algo:fixing} depends on the threshold $\varepsilon_2$: the higher the threshold the larger the number of potentially fixed variables. However, when $\varepsilon_2$ is large, the condition $\alpha_T \geq UB$ becomes harder to be satisfied, and it can be expected to hold only if the current gap at the node is low. We also need to keep into account that we solve problem (\ref{eq:varfix1}) with a first-order method, and hence with a low precision. Indeed, following the approach described in Section \ref{sec:firstorder}, we compare $UB$ with a value that is slightly smaller than $\alpha_T$ since it is derived from a dual feasible solution of problem (\ref{eq:varfix1}).  If $\varepsilon_2$ is too large, the consequence is that no variable is fixed, and the time needed for solving problem (\ref{eq:varfix1}) is wasted. One possibility may be to adapt $\varepsilon_2$ depending on the current gap, making it larger when the gap significantly decreases. However, when the gap is small, in most of the cases the children are pruned, and hence the advantage of the size reduction due to the fixing does not propagate in the tree. In Section \ref{sec:impl_details} we will discuss the choice of $\varepsilon_2$ in our implementation.

All the discussion above refers to multiple variable fixing for variables in $N$ (in fact, all variables in case of binary QPs). For variables in $P$, single variable fixing is possible. For some $i\in P$, if the optimal value of the following SDP problem: 
$$
\begin{array}{ll}
\min_{({\bf X},{\bf x})\in {\cal S}^\star} & {\bf Q}_i{\bf x}+c_i,
\end{array}
$$
is larger than 0, then, according to (\ref{eq:kktboxqp}),
we can fix the value of $x_i$ to 0, while 
if the optimal value of the following SDP problem: 
$$
\begin{array}{ll}
\max_{({\bf X},{\bf x})\in {\cal S}^\star} & {\bf Q}_i{\bf x}+c_i,
\end{array}
$$
is smaller than 0, then, again according to (\ref{eq:kktboxqp}),
we can fix the value of $x_i$ to 1. Note, however, that, due to the large cost of variable fixing, in our experiments we did not include such single variable fixing strategy. 
\begin{algorithm}
\caption{Algorithm for Multiple Variable Fixing}
\label{algo:fixing}
\textbf{MultipleFixing}($UB$, ${\bf x}^\star$, $\varepsilon_1>0$, $\varepsilon_2>0$)
\begin{algorithmic}[1]
\State{{Set $\mathcal{N} = \{x_i : i\in N\}$};}
\While{{$\mathcal{N}$} is not empty}
\State{Set $T_0=\{i\in N\ :\ x_{i}^\star \leq \varepsilon_2\}$ and $T_1=\{i\in N\ :\ x_{i}^\star \geq 1-\varepsilon_2\}$;}
\State{Compute $\alpha_T$ by solving (\ref{eq:varfix1}) with tolerance $\varepsilon_1$;}
\If{$\alpha_T \geq UB$}
\For{$i_r\in \{1, \ldots, \lvert T \rvert\}$}
\If{$i_r\in T_0$}
\State{Fix $x_{i_r}=0$;}
\Else
\State{Fix $x_{i_r}=1$;}
\EndIf
\State{{Set $\mathcal{N} = \mathcal{N} \setminus \{x_{i_r}\}$};}
\EndFor
\State{\textbf{break}}
\Else
\State{{Set $\mathcal{N}=\emptyset$ (variables $x_{i_r}$ for all $i_r \in T$ cannot be fixed);}}
\EndIf
\EndWhile
\end{algorithmic}
\end{algorithm}
\section{Standard features}
\label{sec:standard}
In order to have a full implementation of a branch-and-bound algorithm, we still need to clarify some points, in particular: i) how we initialize and update the upper bound value $UB$; ii) how we perform branching. 
In both cases we use standard rules, which are briefly summarized in this section. Before starting the B\&B algorithm, the upper bound value $UB$ is initialized by running a prefixed number of local searches from randomly
generated points within the unit box. Each time we compute a lower bound, we start a local search from the ${\bf x^\star}$ part of the optimal solution of the relaxation, and we possibly update the upper bound if this is improved by the newly detected local minimizer.
For Binary QP problems, instead, we use the following randomized rounding procedure \cite{goemans1995improved,park2018semidefinite}. From the optimal solution $(\mathbf{X}^\star, \mathbf{x}^\star)$ of the relaxation, we randomly draw several samples $\mathbf{w}$ from $\mathcal{N}(\mathbf{x}^\star, \mathbf{X}^\star - \mathbf{x}^\star (\mathbf{x}^\star)^\top)$, round $w_i$ to 0 or 1 to obtain $\bar{\mathbf{x}} \in \{0, 1\}^n$, and keep the  $\bar{\mathbf{x}}$ with the smallest objective value. This procedure works better than just naively rounding the coordinates of $\mathbf{x}^\star$ to the closest integer.
\noindent
Concerning branching, we select the branching variable according to the following rules. Let $({\bf X}^\star,{\bf x}^\star)$ be the optimal solution of the SDP relaxation at a given node.
First, we search branching variables within the subset $N$ and we select the one furthest from a binary value, i.e., we select $x_i$ such that
\begin{equation*}
i\in\arg\max_{j\in N} \min\{x_j^\star,1-x_j^\star\}.
\end{equation*}
Then, binary branching is performed over the selected variable. If no more variables in $N$ are available for branching, then a variable in $P$ maximizing the approximation error is selected, i.e., 
\begin{equation*}
i\in \arg\max_{j\in P} \sum_{k=1}^n Q_{kj}(x_j^\star x_k^\star - X_{jk}^\star),
\end{equation*}
and ternary branching, as discussed in Section \ref{sec:boxqp}, is performed over the selected variable.
{As already discussed in that section, binary/ternary branching is relevant from a theoretical point of view, guaranteeing equivalence of upper and lower bounds after a finite number of iterations. However, we stress that, from a practical point of view, in our computational experiments over instances taken from the existing literature, B\&B nodes were always fathomed before performing ternary branching on variables in $P$. Thus, to assess the impact of ternary branching, we carried out specific experiments on newly created instances, where $N$ is empty and hence $|P| = n$. The computational results are detailed in Section \ref{sec:ternary_experiments}.}

In conclusion, at any level of the B\&B tree, we obtain lower bounds by solving SDP relaxations of smaller dimension where the size reduction is achieved by fixing decision variables to 0 or 1 thanks to the branching decisions and to the multiple variable fixing strategy.

\section{Computational experiments}
In this section, we describe the implementation details of the B\&B algorithm to globally solve BoxQPs, and we compare its performance with other methods
on test problems taken from the literature. The code has been made publicly available at \url{https://github.com/antoniosudoso/bb-boxqp-fixing}.
\label{sec:exp}
\subsection{Implementation details}
\label{sec:impl_details}
Our B\&B algorithm is implemented in C++
with some routines written in {\tt MATLAB} (version R2020b Update 7).  We use {\tt SDPNAL+} (version 1.0) \cite{sun2020sdpnal+}, a {\tt MATLAB} software that implements an augmented Lagrangian method to solve SDPs with bound constraints. We set the accuracy tolerance of {\tt SDPNAL+} to $10^{-4}$ in the relative KKT residual and we use {\tt Gurobi} (version 10.0.2) to solve the auxiliary LP for producing a valid lower bound. Before starting the B\&B algorithm we perform 200 local searches via {\tt SNOPT} (version 7.7) \cite{gill2005snopt} with randomly generated points within the unit box.
All the experiments are performed on a laptop with an Intel(R) Core(TM) i7-8565U processor clocked at 1.80GHz with 4 cores, 16 GB of RAM and Ubuntu 20.04 LTS. {We provide both a serial implementation, operating on a single thread, and a parallel version with multi-threading capabilities for concurrent node exploration}. 
As for the cutting-plane setting, at each iteration, we separate at most 100000 triangle inequalities, we sort them in decreasing order with respect to the violation and we add the first 10\%. {Therefore, we add at most 10000 violated inequalities at each cutting-plane iteration}.
In our numerical tests, the tolerance for checking the violation
is set to $10^{-4}$. Furthermore, the tolerance for removing inactive triangles is set to $10^{-3}$. {Using this setting, we retain a reasonable number of inequalities at the nodes without compromising the quality of the bound}. We stop the cutting-plane procedure when the number of violated triangles is less than $10n$ or when the relative difference between consecutive lower bounds is less than $10^{-4}$. {We define the relative optimality gap as
\begin{equation}\label{eq:gap}
    \textrm{GAP} = \frac{UB - LB}{|UB|} \times 100,
\end{equation}
where $UB$ is the current global upper bound and $LB$ is the lower bound at the node. A node is fathomed when $\textrm{GAP} \leq 0.01\%$.}
Finally, we explore the B\&B tree with the best-first search strategy.
\subsection{Test instances}
For the computational study, we consider {117} BoxQP instances denoted as {\tt s$\{n\}$-$\{d\}$-$\{s\}$} where ${\bf Q}\in \mathbb{R}^{n\times n}$ is symmetric and not positive semidefinite, ${\bf c}\in \mathbb{R}^n$, $d$ is the density expressed as a percentage, and $s$ is the seed used for the random number generator. These come from the following sources: (i) 36 instances with size $n \in \{70, 80, 90, 100\}$ generated in \cite{Burer09a}, (ii) 9 instances with size $n=125$ generated in \cite{Burer10}; and {(iii)} 27 instances with size $n \in \{150, 175, 200\}$ generated for this paper. The instances with $n \geq 150$ have been generated in the same way as the smallest ones, i.e., for varying densities $d \in \{25, 50, 75\}$ of $\mathbf{Q}$, nonzeros of $\mathbf{Q}$ and $\mathbf{c}$ are uniformly generated integers in $[-50, 50]$. For each pair of $n$ and $d$, three instances have been generated by varying $s \in \{1, 2, 3\}$.
{Furthermore, an additional set of 42 instances is generated, featuring sizes $n \in \{70, 80, 90, 100, 125\}$ and densities $d \in \{25, 50, 75\}$. Notably, for this set, the diagonal elements $Q_{ii}$ are randomly selected as positive integers within the range $[1, 51]$, ensuring $N = \emptyset$ and $|P| = n$. For each pair of $n$ and $d$, three instances are generated, varying the seed $s \in \{4, 5, 6\}$.}

Finally, we consider 40 Binary QP instances available
in the BiqMac library \cite{wiegele2007biq} denoted as {\tt be$\{n\}$-$\{d\}$-$\{s\}$}, with dimension $n \in \{120, 150\}$, density $d \in \{30, 80\}$, and seed $s \in \{1, \dots, 10\}$. As mentioned in Section \ref{sec:boxqp}, we transform Binary QP instances into BoxQP problems by selecting a vector $\boldsymbol{\lambda} \in \mathbb{R}^n$ in such a way that $\lambda_i \leq Q_{ii}$, see Section \ref{sec:binqp} for details.

\subsection{Comparison between Mosek and SDPNAL+}
\label{sec:moseksdpnal}

As a first experiment, we compare, in terms of computing times and quality of the lower bound, the SDP solvers {\tt Mosek} and {\tt SDPNAL+}. We call {\tt Mosek} (version 9.3.20) via {\tt YALMIP} (version 31-March-2021) modeling interface. {Table \ref{tab:my_label4} presents the results at the root node for each size $n\in\{70,80,90,100,125,150\}$, percentage density $d\in\{25,50,75\}$ and seed $s \in \{1, 2, 3\}$. Specifically, we report the following quantities: i) $\textrm{Gap}_0$, the percentage gap according to \eqref{eq:gap}, where LB is computed by solving the root relaxation with only RLTs; ii) CP, the number of cutting-plane iterations where triangle inequalities are added; iii) $\textrm{Gap}_{\textrm{CP}}$, the percentage gap according to \eqref{eq:gap}, where LB is obtained after CP iterations; and iv) Time, computational time in seconds.}
In terms of computing times, {\tt Mosek} and {\tt SDPNAL+} exhibit similar performance up to dimension $n=80$. However, as the dimension increases, {\tt Mosek} becomes significantly slower than {\tt SDPNAL+}. Regarding the quality of the bounds,
we observe that the percentage gaps produced by {\tt Mosek} are slightly lower than those produced by {\tt SDPNAL+}. Indeed, the safe bound of {\tt SDPNAL+}, based on a feasible dual solution, is not necessarily the best one. However, the lower quality of the bound is largely compensated by much lower computing times. The lower quality of the bound also influences the number of cutting-plane iterations, which are generally higher with {\tt SDPNAL+}. Nevertheless, this additional effort is justified by the decreased computational times. 
{To validate this outcome, we have conducted additional testing of our complete B\&B algorithm, employing {\tt Mosek} for the lower bound computation across instances with dimensions up to $n=150$. The results reveal a significantly worse performance in terms of computational time compared to those achieved with {\tt SDPNAL+}. For the sake of space, we do not report these results here. Consequently, for the computational results presented in the upcoming sections, we have chosen to employ the {\tt SDPNAL+} solver and combine it with the post-processing technique detailed in Section \ref{sec:firstorder} to derive valid lower bounds}.
It is also worthwhile to highlight that adding triangle inequalities leads to considerable improvement in the percentage gap at the root node. As could be expected from the improvement of the root lower bound, the cutting-plane procedure with triangle inequalities produces a dramatic improvement in solving BoxQP instances to global optimality.

\begin{table}[!htbp]
    \centering
    \footnotesize
    \caption{{Comparison between Mosek and SDPNAL+ at the root node. We report the following statistics: the percentage gap obtained by solving the basic SDP relaxation ($\textrm{Gap}_0$), the number of cutting-plane iterations (CP), the percentage gap obtained after CP iterations ($\textrm{Gap}_{\textrm{CP}}$), the computational time in seconds (Time).}}
    \label{tab:my_label4}
\begin{tabular*}{\textwidth}{@{\extracolsep{\fill}}cccccccccc@{\extracolsep{\fill}}}
\toprule%
& \multicolumn{4}{c}{Mosek} & \multicolumn{4}{c}{SDPNAL+}\\
\cmidrule{2-5}\cmidrule{6-9}%
Instance & $\textrm{Gap}_0$ & CP & $\textrm{Gap}_{\textrm{CP}}$ & Time &     $\textrm{Gap}_0$ &        CP        &      $\textrm{Gap}_{\textrm{CP}}$   & Time   \\
\midrule
s070-025-1	&	0.23	&	1	&	0.00	&	22	&	0.25	&	1	&	0.00	&	18	\\
s070-025-2	&	1.11	&	1	&	0.00	&	24	&	1.17	&	1	&	0.00	&	19	\\
s070-025-3	&	0.50	&	1	&	0.00	&	22	&	0.52	&	1	&	0.00	&	19	\\
s070-050-1	&	0.79	&	1	&	0.00	&	23	&	0.82	&	1	&	0.00	&	17	\\
s070-050-2	&	0.13	&	1	&	0.00	&	22	&	0.15	&	1	&	0.00	&	17	\\
s070-050-3	&	0.00	&	0	&	0.00	&	8	&	0.00	&	0	&	0.00	&	5	\\
s070-075-1	&	0.32	&	1	&	0.00	&	21	&	0.34	&	1	&	0.00	&	17	\\
s070-075-2	&	2.09	&	1	&	0.00	&	22	&	2.11	&	2	&	0.00	&	40	\\
s070-075-3	&	1.13	&	1	&	0.00	&	22	&	1.15	&	1	&	0.01	&	19	\\
\midrule
s080-025-1	&	0.00	&	0	&	0.00	&	15	&	0.02	&	1	&	0.00	&	13	\\
s080-025-2	&	1.39	&	1	&	0.00	&	28	&	1.41	&	1	&	0.00	&	17	\\
s080-025-3	&	0.42	&	1	&	0.00	&	28	&	0.45	&	1	&	0.00	&	22	\\
s080-050-1	&	3.95	&	2	&	0.11	&	42	&	4.01	&	2	&	0.13	&	58	\\
s080-050-2	&	0.04	&	1	&	0.00	&	27	&	0.05	&	1	&	0.00	&	26	\\
s080-050-3	&	0.71	&	1	&	0.00	&	33	&	0.72	&	1	&	0.00	&	32	\\
s080-075-1	&	0.47	&	1	&	0.00	&	28	&	0.50	&	1	&	0.00 &	24	\\
s080-075-2	&	0.90	&	1	&	0.00	&	30	&	0.91	&	2	&	0.00	&	58	\\
s080-075-3	&	0.97	&	1	&	0.00	&	34	&	1.02	&	2	&	0.00	&	46	\\
\midrule
s090-025-1	&	1.80	&	1	&	0.00	&	42	&	1.83	&	1	&	0.00	&	28	\\
s090-025-2	&	1.40	&	1	&	0.00	&	50	&	1.42	&	1	&	0.00 &	36	\\
s090-025-3	&	0.99	&	1	&	0.00	&	58	&	1.02	&	1	&	0.00	&	35	\\
s090-050-1	&	1.39	&	1	&	0.00	&	61	&	1.45	&	1	&	0.00	&	36	\\
s090-050-2	&	0.00	&	0	&	0.00	&	33	&	0.02	&	1	&	0.00	&	26	\\
s090-050-3	&	0.67	&	1	&	0.00	&	59	&	0.68	&	1	&	0.00	&	39	\\
s090-075-1	&	2.17	&	2	&	0.00	&	89	&	2.19	&	2	&	0.02	&	73	\\
s090-075-2	&	2.32	&	2	&	0.01	&	92	&	2.34	&	2	&	0.07	&	74	\\
s090-075-3	&	1.26	&	1	&	0.00	&	59	&	1.30	&	1	&	0.00	&	30	\\
\midrule
s100-025-1	&	0.97	&	1	&	0.00	&	83	&	0.99	&	1	&	0.00	&	31	\\
s100-025-2	&	0.81	&	1	&	0.00	&	78	&	0.82	&	1	&	0.00	&	35	\\
s100-025-3	&	0.53	&	1	&	0.00	&	81	&	0.54	&	1	&	0.00	&	30	\\
s100-050-1	&	3.31	&	2	&	0.25	&	129	&	3.35	&	2	&	0.28	&	73	\\
s100-050-2	&	2.20	&	2	&	0.08	&	141	&	2.21	&	3	&	0.10	&	84	\\
s100-050-3	&	0.86	&	1	&	0.00	&	86	&	0.88	&	1	&	0.00	&	33	\\
s100-075-1	&	1.76	&	2	&	0.02	&	136	&	1.78	&	3	&	0.04	&	104	\\
s100-075-2	&	1.90	&	2	&	0.09	&	129	&	1.96	&	2	&	0.11	&	74	\\
s100-075-3	&	1.69	&	2	&	0.00	&	124	&	1.71	&	3	&	0.01	&	100	\\
\midrule
s125-025-1	&	2.84	&	2	&	0.40	&	313	&	2.86	&	2	&	0.43	&	102	\\
s125-025-2	&	0.49	&	1	&	0.00	&	186	&	0.51	&	1	&	0.00	&	41	\\
s125-025-3	&	0.37	&	1	&	0.00	&	208	&	0.38	&	1	&	0.00	&	64	\\
s125-050-1	&	1.95	&	2	&	0.11	&	332	&	1.98	&	3	&	0.14	&	164	\\
s125-050-2	&	0.61	&	1	&	0.00	&	180	&	0.63	&	1	&	0.00	&	50	\\
s125-050-3	&	1.02	&	2	&	0.00	&	318	&	1.05	&	3	&	0.01	&	228	\\
s125-075-1	&	3.41	&	2	&	0.71	&	287	&	3.48	&	2	&	0.74	&	101	\\
s125-075-2	&	2.07	&	2	&	0.26	&	298	&	2.10	&	3	&	0.27	&	157	\\
s125-075-3	&	2.55	&	2	&	0.58	&	282	&	2.57	&	3	&	0.60	&	107	\\
\midrule
s150-025-1	&	4.03	&	2	&	1.48	&	669	&	4.04	&	3	&	1.53	&	134	\\
s150-025-2	&	1.91	&	3	&	0.07	&	1000	&	1.94	&	3	&	0.10	&	183	\\
s150-025-3	&	1.06	&	1	&	0.00	&	533	&	1.07	&	2	&	0.00	&	141	\\
s150-050-1	&	2.17	&	2	&	0.03	&	696	&	2.21	&	3	&	0.06	&	189	\\
s150-050-2	&	0.92	&	2	&	0.0	&	812	&	0.95	&	2	&	0.00	&	173	\\
s150-050-3	&	1.40	&	2	&	0.0	&	764	&	1.43	&	3	&	0.01	&	234	\\
s150-075-1	&	4.43	&	2	&	2.01	&	726	&	4.45	&	2	&	2.08	&	138	\\
s150-075-2	&	1.33	&	3	&	0.01	&	982	&	1.37	&	3	&	0.04	&	207	\\
s150-075-3	&	2.28	&	2	&	0.57	&	714	&	2.31	&	3	&	0.59	&	201	\\
\bottomrule
\end{tabular*}
\end{table}

\subsection{Comparison with existing solvers}
\begin{figure}[!htbp]
    \centering
    \includegraphics[scale=0.7]{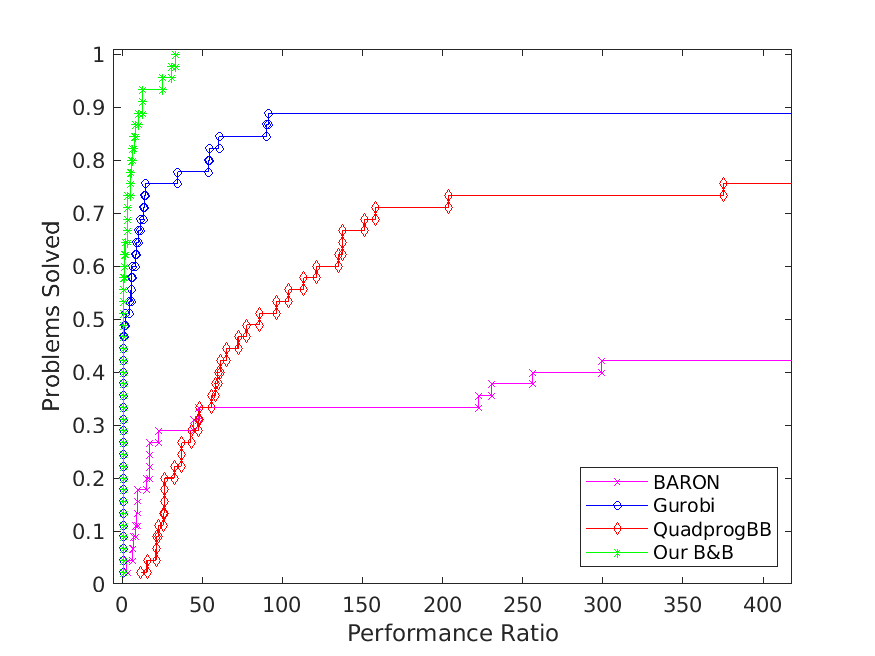}
    \caption{{Computational time performance profile comparing solvers on BoxQP instances ({\tt s$\{n\}$-$\{d\}$-$\{s\}$}, for $n \in \{70, 80, 90, 100, 125\}$, $d \in \{25, 50, 75\}$, and $s \in \{1, 2, 3\}$).}}
    \label{fig:perf_boxQP}
\end{figure}

In this section, we report experimental results designed to compare the performance of our solver
against other state-of-the-art methods for solving non-convex BoxQPs to global optimality.
According to some recent results presented in \cite{Bonami18,Xia20}, the best solvers for BoxQP problems are
{\tt Cplex} and, for instances with a large density, the SDP-based branch-and-bound approach {\tt QuadprogBB} \cite{Chen12}.
In the experiments presented in those works, {\tt Gurobi} was not included. However, according to our experiments, 
the latest version of {\tt Gurobi} performs better than {\tt Cplex}. For this reason, we compared our B\&B approach  
with {\tt Gurobi} {(version 10.0.2)} and {\tt QuadprogBB} {over the {\tt s$\{n\}$-$\{d\}$-$\{s\}$} instances with $s \in \{1, 2, 3\}$ and $n \leq 125$}.
{We compare with {\tt QuadprogBB} since it is the best-known solver for BoxQP problems based on the computation of SDP bounds and, as already mentioned, it performs well over large density instances. 
Moreover, we compare against \texttt{BARON} \cite{sahinidis1996baron,khajavirad2018hybrid}, a general-purpose software that is renowned in the literature for solving non-convex optimization problems to global optimality. We call \texttt{BARON} (version
23.3.11) via AMPL modeling interface. Finally, to ensure a fair and unbiased comparison, we consider the serial (single-threaded) versions of all solvers, including the B\&B proposed in this paper}. Each instance in the test suite is solved with a time limit of 7200 seconds and to a final relative optimality gap tolerance of 0.01\%.
The results are displayed in Table \ref{tab:comparison_spar123}. Whenever the time limit of 2 hours is reached, we report in column ``Nodes'' the final percentage gap when the algorithm stops. Overall, we can make the following observations:
\begin{itemize}
\item {{\tt BARON} and {\tt Gurobi} are extremely efficient over small ($n \in \{70, 80, 90\}$) and low-density $(d=25)$ instances. Here, the high quality of the SDP bounds used by {\tt QuadprogBB} and our solver is not able to compensate for the high computational cost.}
\item {\tt QuadprogBB} outperforms {{\tt BARON}, and becomes more competitive at large densities, as already observed in the literature. It sometimes outperforms {\tt Gurobi} over large density instances. For example, for $n=125$ and $d=75$, the average final optimality gap after two hours is always lower than the one obtained by {\tt Gurobi}.}
\item {{\tt BARON} struggles more than {\tt Gurobi} when both the size and the density increase, reaching the time limit on all the instances with $n=125$.}
\item Our B\&B approach always outperforms {\tt QuadprogBB}, it strongly outperforms {\tt Gurobi} {and {\tt BARON}} at density $d=75$, and the difference becomes more and more evident as the dimension increases; at density $d=50$, {\tt Gurobi} performs better at low dimension but our B\&B becomes better at higher dimensions; at density $d=25$, {\tt Gurobi} is almost always better with the remarkable exception of the instance {\tt s125-025-1}.
\item In our approach, the number of nodes is much lower than the number of nodes generated by the other {three} approaches and, for all the instances with the exception of {\tt s125-050-1}, the globally optimal solution is already found (and sometimes also certified as such) at the root node.
\end{itemize}
{For the sake of fairness, we point out that the poor performance of {\tt BARON} on dense instances, when compared to other solvers, is not surprising. Indeed, {\tt BARON} is a general-purpose global optimization solver, while all the other solvers are tailored specifically for QP problems. Moreover, we point out that}
{\tt QuadprogBB} has not been further developed in recent years. Therefore, we believe that
its computing times can be improved. For instance, it could be enhanced through the use of 
a first-order solver like {\tt SDPNAL+}. But we remark that, in terms of computing times, our approach 
outperforms {\tt QuadprogBB} even when we employ the SDP solver {\tt Mosek}. For example, 
all instances {\tt s125-075}, which are not solved within 7200 seconds by {\tt QuadprogBB}, are solved in at most 1642 {seconds} by our B\&B approach without multiple variable fixing, and in at most 994 seconds by our approach with multiple variable fixing, when SDP subproblems are solved with {\tt Mosek}. Therefore, it seems 
that what really makes the difference is the intensive introduction of cutting planes and the newly proposed multiple fixing procedure. 
Concerning the latter, we point out that in Table \ref{tab:comparison_spar123} we have not reported results with such procedure. On instances with size up to $n=100$, where a very low number of nodes is generated, the impact of multiple fixing is mild. But as we will see later on, this impact becomes quite significant at larger dimensions. In Figure \ref{fig:perf_boxQP} we report a performance profile \cite{dolan2002benchmarking} of the computational time comparing the solvers over all the {{\tt s$\{n\}$-$\{d\}$-$\{s\}$} instances with $s \in \{1, 2, 3\}$ and $n \leq 125$, which clearly shows the effectiveness of our B\&B algorithm}.

{As a final remark about these tests, we observe that, while we compare against some well-known publicly available solvers, the literature suggests other approaches for BoxQP problems. For instance, in \cite{elloumi2019global}, a method for addressing mixed-integer QCQP problems is presented. This study also includes tests on BoxQP problems, showing slightly better results than those obtained with \texttt{BARON}.}

\begin{table}[!htbp]
\centering
\caption{{Comparison between BARON, Gurobi, QuadprogBB and our B\&B solver on BoxQP instances. For each instance, we report the number of nodes (Nodes) and the computational time in seconds (Time). When the time limit of 7200 seconds is reached we report the final percentage gap in the column ``Nodes''.}} 
\footnotesize
\begin{tabular}{lrr|rr|rr|rr}
\toprule%
& \multicolumn{2}{@{}c@{}}{BARON} & \multicolumn{2}{@{}c@{}}{Gurobi} & \multicolumn{2}{@{}c@{}}{{QuadprogBB}} & \multicolumn{2}{@{}c@{}}{Our B\&B}\\
\cmidrule{2-3}\cmidrule{4-5}\cmidrule{6-7}\cmidrule{8-9}%
Instance &  Nodes & Time & Nodes &  Time &  Nodes &  Time & Nodes & Time\\
\midrule
s070-025-1 &      5 &     5 &         13 &      1 &     15 &   56 &      1 &    18 \\
s070-025-2 &      4 &     5 &         28 &      1 &     17 &   90 &      1 &    20 \\
s070-025-3 &      5 &     5 &         10 &      1 &     19 &  119 &      1 &    18 \\
s070-050-1 &     14 &   208 &         18 &      5 &     19 &  105 &      1 &    16 \\
s070-050-2 &      6 &    21 &          3 &      2 &     17 &  102 &      1 &    17 \\
s070-050-3 &      1 &     3 &          1 &      1 &      1 &   25 &      1 &     8 \\
s070-075-1 &    (5.14\%) &  7200 &         12 &     29 &     23 &  187 &      1 &    16 \\
s070-075-2 &    (5.30\%) &  7200 &        520 &    206 &    149 &  875 &      1 &    33 \\
s070-075-3 &    (6.64\%) &  7200 &        162 &    109 &     63 &  410 &      1 &    19 \\
\midrule
s080-025-1 &      1 &     1 &          2 &      1 &     11 &   97 &      1 &    13 \\
s080-025-2 &     30 &   238 &        389 &      5 &     33 &  278 &      1 &    17 \\
s080-025-3 &      5 &    14 &         19 &      2 &     29 &  286 &      1 &    20 \\
s080-050-1 &    (6.46\%) &  7200 &       1050 &    350 &    455 & 3713 &      3 &    62 \\
s080-050-2 &      4 &    96 &          1 &      6 &      7 &   98 &      1 &    21 \\
s080-050-3 &      3 &   138 &         22 &     14 &     51 &  527 &      1 &    27 \\
s080-075-1 &    (29.53\%) &  7200 &         88 &    287 &     83 &  786 &      1 &    21 \\
s080-075-2 &    (24.26\%) &  7200 &        200 &    207 &     99 &  914 &      1 &    43 \\
s080-075-3 &    (20.63\%) &  7200 &        551 &    432 &    209 & 2578 &      1 &    42 \\
\midrule
s090-025-1 &      6 &    77 &        216 &      4 &     77 &  709 &      1 &    24 \\
s090-025-2 &     10 &    69 &        244 &      4 &     41 &  550 &      1 &    27 \\
s090-025-3 &     10 &    49 &         34 &      2 &     87 &  808 &      1 &    27 \\
s090-050-1 &    (4.67\%) &  7200 &        108 &     74 &    195 & 2389 &      1 &    33 \\
s090-050-2 &     23 &   270 &          1 &     15 &     29 &  413 &      1 &    21 \\
s090-050-3 &     28 &  7099 &         29 &     24 &     83 & 1133 &      1 &    32 \\
s090-075-1 &    (82.01\%) &  7200 &       1019 &   4117 &    453 & 5285 &      3 &    68 \\
s090-075-2 &    (86.48\%) &  7200 &        925 &   3713 &    403 & 4439 &      3 &    68 \\
s090-075-3 &    (64.76\%) &  7200 &        693 &   2467 &    135 & 1573 &      1 &    27 \\
\midrule
s100-025-1 &     61 &  6050 &        680 &     24 &     45 &  621 &      1 &    28 \\
s100-025-2 &     19 &  1916 &        234 &      9 &     23 &  371 &      1 &    28 \\
s100-025-3 &     19 &  1087 &         85 &      5 &     41 &  649 &      1 &    27 \\
s100-050-1 &    (14.64\%) &  7200 &       3675 &   2610 &    (1.16\%) & 7200 &      5 &   186 \\
s100-050-2 &    (7.85\%) &  7200 &       1257 &    925 &    (0.17\%) & 7200 &      3 &   107 \\
s100-050-3 &    (7.88\%) &  7200 &        251 &    275 &     67 &  976 &      1 &    30 \\
s100-075-1 &    (111.55\%) &  7200 &       (3.19\%) &   7200 &    (0.13\%) & 7200 &      3 &   110 \\
s100-075-2 &    (111.57\%) &  7200 &       (2.24\%) &   7200 &    (0.46\%) & 7200 &      3 &   122 \\
s100-075-3 &    (98.91\%) &  7200 &        908 &   6664 &    (0.31\%) & 7200 &      1 &    74 \\
\midrule
s125-025-1 &   (6.94\%) &  7200 &      36302 &   4007 &    (11.81\%) & 7200 &      7 &   340 \\
s125-025-2 &     (3.93\%) &  7200 &        548 &     52 &    213 & 5097 &      1 &    42 \\
s125-025-3 &     (1.31\%) &  7200 &         48 &     11 &     69 & 1485 &      1 &    58 \\
s125-050-1 &    (69.38\%) &  7200 &       1745 &   6476 &    (7.59\%) & 7200 &      3 &   186 \\
s125-050-2 &    (56.87\%) &  7200 &        715 &   2536 &    137 & 4025 &      1 &    47 \\
s125-050-3 &    (53.04\%) &  7200 &        605 &   2511 &     (0.15\%) & 7200 &      1 &   171 \\
s125-075-1 &    (200.00\%) &  7200 &        (25.96\%) &   7200 &    (17.54\%) & 7200 &     13 &   489 \\
s125-075-2 &    (163.57\%) &  7200 &        (15.98\%) &   7200 &    (11.46\%) & 7200 &      7 &   435 \\
s125-075-3 &    (158.75\%) &  7200 &        (16.97\%) &   7200 &    (14.23\%) & 7200 &     15 &   719 \\
\bottomrule
\end{tabular}
    \label{tab:comparison_spar123}
\end{table}

\subsection{The impact of multiple variable fixing}
In this section, we report the results of an experiment designed to evaluate the performance of our B\&B approach with the introduction of multiple variable fixing. 
{The experimental setup for this study is the one described in Section \ref{sec:impl_details}. However, due to the large size of the considered instances, we use the parallel implementation of our B\&B method to speed up the computational times. In this setting, we process two child nodes in parallel by assigning two threads to each of them. }
For the multiple variable fixing, we follow Algorithm \ref{algo:fixing} taking as input the $\mathbf{x}$ part of the solution of the SDP relaxation solved at each node, the current $UB$, and $\varepsilon_1 =\varepsilon_2=0.01$. We only activate the fixing procedure when the gap at the node is less than 1\%. This strategy turns out to be effective at every node since no time is wasted for fixing variables that cannot be fixed since the node gap is too large.

{Table \ref{tab:comparison_fixing} compares our B\&B approach without and with multiple fixing over 36 {\tt s$\{n\}$-$\{d\}$-$\{s\}$} instances where $n \in \{125,150,175,200\}$, $d \in \{25, 50, 75\}$, and $s \in \{1, 2, 3\}$. The first column, $\vert N \vert$, represents the cardinality of set $N$; columns ``Nodes'' and ``Time'' indicate the number of nodes and computational time in seconds, respectively; column ``Node Red.'' expresses the percentage reduction in nodes for B\&B with fixing compared to B\&B without fixing. This reduction is calculated as 
$100\frac{C_1-C_2}{C_1}$, where $C_1$ denotes the number of nodes without fixing and $C_2$ the number of nodes with fixing; similarly, column ``Time Red.'' denotes the percentage reduction in computational time for B\&B with fixing compared to B\&B without fixing, calculated as 
$100\frac{T_1-T_2}{T_1}$, where $T_1$ denotes 
the time without fixing and $T_2$ the time with fixing.}

These computational results clearly show that multiple fixing, despite its additional cost (each round of fixing requires the solution of an SDP problem), is highly beneficial, especially over the most challenging instances, when the number of nodes becomes larger. As already commented in Section \ref{sec:fixing}, the reduction of the computing times is partially due to the usually lower number of nodes, but is mostly related to the smaller dimension of the SDP problems at the B\&B nodes{, as confirmed by the statistics in Table \ref{tab:fix_stat}.} This is clearly highlighted by instance {\tt s200-050-1}, {where the number of nodes with multiple variable fixing increases of $8.36\%$ (this happens because the choice of branching variables is different in the two approaches), but the computational time reduces of $32.97\%$. By comparing the two columns ``Node Red.'' and ``Time Red.'', it becomes clear that the reduction of computing times cannot only be explained by the reduction of the number of nodes. Remarkably, on the largest instances ($n=200$) we can observe a time reduction greater than $30\%$ with the exception of \texttt{s200-025-3} which is solved at the root node without fixing any variable.}

{To better analyze the impact of the fixing procedure, we report in Table \ref{tab:fix_stat} some statistics. These include the total number of SDPs solved for fixing variables (\# SDP); the computational time in seconds required for solving these SDPs (Time); the number of ``wasted'' SDPs (\# W), that is the SDPs that did not allow to fix any variable; the total time in seconds spent on the wasted SDPs (W Time); the number of variables fixed at the root node (Fix$_0$); the average number of fixed variables (Avg); the minimum (Min) and maximum (Max) number of fixed variables. Finally, we report the average size $n$ of the SDPs solved at the nodes (Avg $n$).   
Looking at Table \ref{tab:fix_stat}, it emerges that the number of wasted SDPs is in general neglectable with respect to the total number of solved SDPs. Even in scenarios where solving some SDPs does not allow to fix variables, such as in \texttt{s200-050-2}, where fixing fails in 10 out of 63 cases, the resulting average size of 171.53 leads to an impressive time reduction of $43.83\%$ (see Table \ref{tab:comparison_fixing}). We highlight that the high quality of the SDP bound is crucial for the success of the multiple variable fixing. Indeed, we try to fix some variables at each node only when the relative gap is below $1\%$. Column ``Avg $n$'' in Table \ref{tab:fix_stat} shows a significant drop in the size at the nodes, which implies that the gap goes below $1\%$ early enough in the tree, thanks to the quality of the bound. Once again, the column ``Avg $n$'' confirms that the primary impact of multiple fixing lies in the reduction of SDP sizes solved within the branch-and-bound tree. This size reduction is particularly beneficial given the sensitivity of state-of-the-art SDP solvers to the problem dimension.}

\begin{table}[!htbp]
    \centering
    \footnotesize
    \caption{{Comparison between our B\&B solver without and with multiple variable fixing. Column $\vert N \vert$ denotes the cardinality of set $N$; columns ``Nodes'' and ``Time'' indicate the number of nodes and computational time in seconds, respectively; column ``Node Red.'' expresses the percentage reduction in nodes for B\&B with fixing compared to B\&B without fixing; column ``Time Red.'' denotes the percentage reduction in computational time for B\&B with fixing compared to B\&B without fixing.}}
\begin{tabular}{lccc|cccc}
\toprule%
\multicolumn{2}{c}{} & \multicolumn{2}{c}{Our B\&B} & \multicolumn{4}{c}{Our B\&B with Fixing}\\
\cmidrule{3-4}\cmidrule{5-8}%
Instance	& $\vert N \vert$ &		Nodes	&	Time	&	Nodes	&	Time &	Node Red. [\%]	&	Time Red. [\%]	\\
\midrule
s125-025-1	&	105 &	7	&	309	&	7	&	262	&	0	&	15.21	\\
s125-025-2	&	109 &	1	&	41	&	1	&	41	&	0	&	0	\\
s125-025-3	&	113 &	1	&	63	&	1	&	61	&	0	&	3.17	\\
s125-050-1	&	97 &	3	&	222	&	3	&	224	&	0	&	-0.9	\\
s125-050-2	&	100 &	1	&	46	&	1	&	45	&	0	&	2.17	\\
s125-050-3	&	95 &	1	&	190	&	1	&	192	&	0	&	-1.05	\\
s125-075-1	&	75 &	13	&	409	&	13	&	360	&	0	&	11.98	\\
s125-075-2	&	74 &	7	&	404	&	7	&	313	&	0	&	22.52	\\
s125-075-3	&	74 &	15	&	657	&	15	&	397	&	0	&	39.57	\\
\midrule
s150-025-1	&	127 &	85	&	4249	&	73	&	2511	&	14.12	&	40.9	\\
s150-025-2	&	132 &	3	&	291	&	3	&	273	&	0	&	6.19	\\
s150-025-3	&	135 &	1	&	138	&	1	&	137	&	0	&	0.72	\\
s150-050-1	&	116 &	3	&	228	&	3	&	237	&	0	&	-3.95	\\
s150-050-2	&	119 &	1	&	139	&	1	&	138	&	0	&	0.72	\\
s150-050-3	&	115 &	1	&	198	&	1	&	194	&	0	&	2.02	\\
s150-075-1	&	91 &	161	&	6742	&	139	&	4266	&	13.66	&	36.73	\\
s150-075-2	&	92 &	3	&	265	&	3	&	269	&	0	&	-1.51	\\
s150-075-3	&	88 &	21	&	1376	&	21	&	807	&	0	&	41.35	\\
\midrule
s175-025-1	&	151 &	11	&	779	&	9	&	632	&	18.18	&	18.87	\\
s175-025-2	&	151 &	7	&	789	&	7	&	791	&	0	&	-0.25	\\
s175-025-3	&	156 &	1	&	207	&	1	&	209	&	0	&	-0.97	\\
s175-050-1	&	137 &	29	&	1810	&	23	&	1073	&	20.69	&	40.72	\\
s175-050-2	&	139 &	9	&	899	&	7	&	470	&	22.22	&	47.72	\\
s175-050-3	&	134 &	1	&	179	&	1	&	181	&	0	&	-1.12	\\
s175-075-1	&	103 &	921	&	47995	&	867	&	32888	&	5.86	&	31.48	\\
s175-075-2	&	109 &	3	&	451	&	3	&	333	&	0	&	26.16	\\
s175-075-3	&	106 &	9	&	773	&	7	&	554	&	22.22	&	28.33	\\
\midrule
s200-025-1	&	175 &	37	&	3471	&	35	&	2338	&	5.41	&	32.64	\\
s200-025-2	&	173 &	31	&	2988	&	27	&	1870	&	12.9	&	37.42	\\
s200-025-3	&	179 &	1	&	354	&	1	&	355	&	0	&	-0.28	\\
s200-050-1	&	158 &	359	&	28134	&	389	&	18858	&	-8.36	&	32.97	\\
s200-050-2	&	156 &	151	&	13478	&	131	&	7571	&	13.25	&	43.83	\\
s200-050-3	&	154 &	101	&	10583	&	95	&	5875	&	5.94	&	44.49	\\
s200-075-1	&	123 &	833	&	71334	&	773	&	40890	&	7.2	&	42.68	\\
s200-075-2	&	126 &	17	&	1850	&	15	&	1145	&	11.76	&	38.11	\\
s200-075-3	&	123 &	123	&	9508	&	105	&	6424	&	14.63	&	32.44	\\
\bottomrule
\end{tabular}
    \label{tab:comparison_fixing}
\end{table}

\begin{table}[!ht]
    \centering
    \footnotesize
    \caption{{Statistics our B\&B with multiple variable fixing. We report the instance name, the number of SDPs solved for fixing (\# SDP), the time in seconds for solving the SDPs for fixing, the number of ``wasted'' SDPs (\# W), that is the number of times where SDPs failed to fix any variable, the time in seconds spent for solving the wasted SDPs (W Time), the number of variables fixed at the root node (Fix$_0$), the average number of fixed variables (Avg), the minimum number of fixed variables (Min) and the maximum number of fixed variables (Max), and the average dimension $n$ of the SDPs solved at the nodes (Avg $n$).}}
\begin{tabular*}{\textwidth}{@{\extracolsep{\fill}}lccccccccc@{\extracolsep{\fill}}}
\toprule%
Instance	&	\# SDP	&	Time	&	\# W	&	W Time	&	Fix$_0$	&	Avg	&	Min	&	Max	&	Avg $n$	\\\midrule
s125-025-1	&	3	&	74	&	0	&	0	&	16	&	9.67	&	3	&	16	&	106.14	\\
s125-025-2	&	0	&	0	&	0	&	0	&	0	&	0.0	&	0	&	0	&	125.0	\\
s125-025-3	&	0	&	0	&	0	&	0	&	0	&	0.0	&	0	&	0	&	125.0	\\
s125-050-1	&	1	&	33	&	0	&	0	&	40	&	40.0	&	40	&	40	&	97.67	\\
s125-050-2	&	0	&	0	&	0	&	0	&	0	&	0.0	&	0	&	0	&	125.0	\\
s125-050-3	&	0	&	0	&	0	&	0	&	0	&	0.0	&	0	&	0	&	125.0	\\
s125-075-1	&	4	&	66	&	0	&	0	&	9	&	5.0	&	2	&	9	&	112.08	\\
s125-075-2	&	3	&	80	&	0	&	0	&	20	&	14.67	&	7	&	20	&	99.57	\\
s125-075-3	&	5	&	89	&	0	&	0	&	12	&	4.6	&	2	&	12	&	110.07	\\
\midrule
s150-025-1	&	33	&	994	&	4	&	147	&	0	&	7.45	&	0	&	20	&	129.51	\\
s150-025-2	&	1	&	56	&	0	&	0	&	46	&	46.0	&	46	&	46	&	118.67	\\
s150-025-3	&	0	&	0	&	0	&	0	&	0	&	0.0	&	0	&	0	&	150.0	\\
s150-050-1	&	1	&	41	&	0	&	0	&	26	&	26.0	&	26	&	26	&	132.0	\\
s150-050-2	&	0	&	0	&	0	&	0	&	0	&	0.0	&	0	&	0	&	150.0	\\
s150-050-3	&	0	&	0	&	0	&	0	&	0	&	0.0	&	0	&	0	&	150.0	\\
s150-075-1	&	44	&	1285	&	0	&	0	&	0	&	4.45	&	2	&	13	&	135.21	\\
s150-075-2	&	1	&	45	&	0	&	0	&	32	&	32.0	&	32	&	32	&	128.0	\\
s150-075-3	&	9	&	260	&	0	&	0	&	13	&	9.67	&	2	&	39	&	125.71	\\
\midrule
s175-025-1	&	5	&	210	&	0	&	0	&	16	&	13.8	&	3	&	31	&	153.44	\\
s175-025-2	&	3	&	260	&	1	&	61	&	0	&	33.33	&	0	&	51	&	145.0	\\
s175-025-3	&	0	&	0	&	0	&	0	&	0	&	0.0	&	0	&	0	&	175.0	\\
s175-050-1	&	10	&	375	&	0	&	0	&	14	&	10.6	&	2	&	51	&	146.65	\\
s175-050-2	&	3	&	96	&	0	&	0	&	30	&	18.33	&	7	&	30	&	140.71	\\
s175-050-3	&	0	&	0	&	0	&	0	&	0	&	0.0	&	0	&	0	&	175.0	\\
s175-075-1	&	253	&	8985	&	0	&	0	&	0	&	3.36	&	2	&	11	&	157.41	\\
s175-075-2	&	1	&	52	&	0	&	0	&	32	&	32.0	&	32	&	32	&	153.0	\\
s175-075-3	&	3	&	121	&	0	&	0	&	15	&	9.67	&	2	&	15	&	157.0	\\
\midrule
s200-025-1	&	16	&	807	&	2	&	110	&	0	&	8.56	&	0	&	21	&	176.23	\\
s200-025-2	&	13	&	675	&	2	&	101	&	0	&	10.92	&	0	&	36	&	171.85	\\
s200-025-3	&	0	&	0	&	0	&	0	&	0	&	0.0	&	0	&	0	&	200.0	\\
s200-050-1	&	128	&	6929	&	2	&	96	&	0	&	3.84	&	0	&	17	&	178.93	\\
s200-050-2	&	63	&	2730	&	10	&	560	&	0	&	8.29	&	0	&	28	&	171.53	\\
s200-050-3	&	45	&	2111	&	6	&	392	&	0	&	8.89	&	0	&	30	&	170.11	\\
s200-075-1	&	273	&	13643	&	6	&	365	&	0	&	3.29	&	0	&	14	&	179.02	\\
s200-075-2	&	6	&	274	&	0	&	0	&	19	&	9.67	&	2	&	19	&	172.13	\\
s200-075-3	&	49	&	2405	&	0	&	0	&	0	&	4.04	&	2	&	23	&	182.53	\\
\bottomrule
\end{tabular*}
    \label{tab:fix_stat}
\end{table}

\subsection{Ternary branching}
\label{sec:ternary_experiments}
{
When evaluating our B\&B algorithm on BoxQP instances, nodes are always fathomed before performing ternary branching on variables within set $P$. Therefore, in this section, to evaluate the effectiveness of ternary branching, we perform some experiments on instances with $N=\emptyset$ ($|P| = n$). 
As before, we compare the serial version of our B\&B solver with the serial versions of {\tt BARON}, {\tt Gurobi}, and {\tt QuadprogBB} over the ${\tt s}\{n\}-\{d\}-\{s\}$ instances with seed $s\in \{4, 5, 6\}$, dimension $n\in\{70,80,90,100,125\}$ and density $d\in \{25,50,75\}$. Each instance is solved with a time limit of 7200 seconds and to a final relative optimality gap tolerance of 0.01\%. The results are shown in Table \ref{tab:ternary_results}. Whenever the time limit of 2 hours is reached, we report in column ``Nodes'' the final percentage gap when the algorithm stops. It turns out that our algorithm still outperforms the other solvers. It is interesting to remark that {\tt Gurobi} performs worse over these instances compared to the previous ones (where $N$ is not empty), while {\tt QuadprogBB} is often better than {\tt Gurobi} at larger dimensions. {\tt BARON} does not perform well, but we emphasize once again that it serves as a general-purpose global optimization solver, unlike the others that are designed specifically for QP problems. As for our B\&B algorithm, we successfully solve all instances within the specified time limit, and the majority of them are solved at the root node. However, the computing times of the hardest instances are significantly larger than those of the hardest instances displayed in Table \ref{tab:comparison_spar123}. We remark that, since $|N|=\emptyset$, for these instances we could not adopt the fixing procedure. However, the fixing procedure could be employed at some nodes. Indeed, for each $i$, let us denote by
\begin{align*}
    N(i)=\{j\neq i\ :\ Q_{ii}Q_{jj}-Q_{ij}^2 \leq 0\}.
\end{align*}
For each $j\in N(i)$, at optimal solutions of the BoxQP problem, either $x_i$ or $x_j$ have a binary value. Therefore, if at some node of the B\&B tree we set ${\bf Q}_i{\bf x}=-c_i$ (i.e., we impose $0<x_i<1$), then at such node we can impose the binary condition for all variables in $N(i)$. However, we have not explored this possibility in the current experiments.}

\begin{figure}[!htbp]
    \centering
    \includegraphics[scale=0.7]{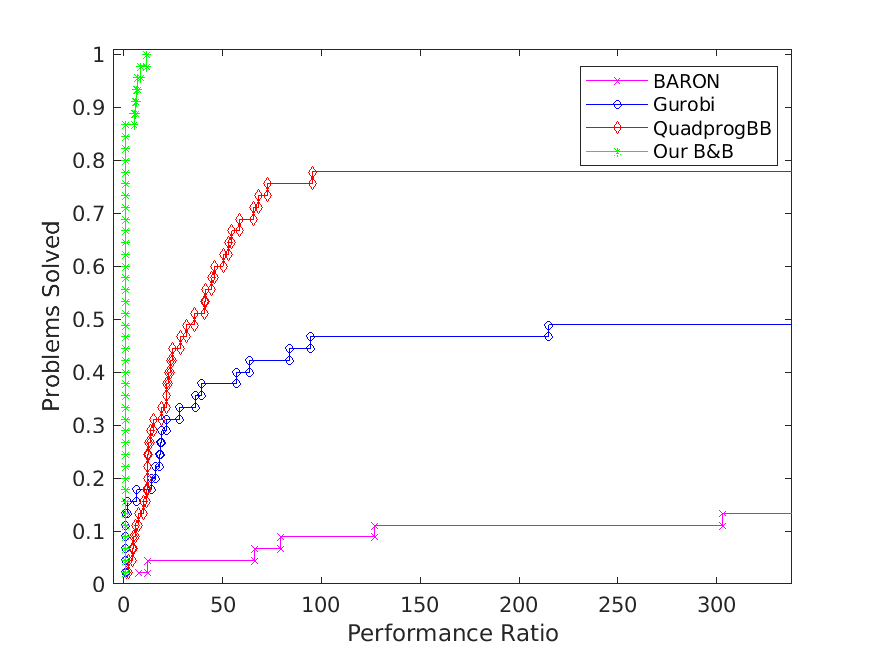}
    \caption{{Computational time performance profile comparing solvers on the BoxQP instances where $N=\emptyset$ ({\tt s$n$-$d$-$s$}, for $n \in \{70, 80, 90, 100, 125\}$, $d \in \{25, 50, 75\}$, and $s \in \{4, 5, 6\}$).}}
    \label{fig:perf_boxQP_emptyN}
\end{figure}

{In Figure \ref{fig:perf_boxQP_emptyN} we report the performance profile of the computational time comparing the solvers over all the {{\tt s$\{n\}$-$\{d\}$-$\{s\}$} instances with $n \leq 125$, $d \in \{25, 50, 75\}$, and $s \in \{4, 5, 6\}$, which clearly shows the effectiveness of our B\&B algorithm}.}

\begin{table}[!htbp]
\centering
\caption{{Comparison between BARON, Gurobi, QuadprogBB and our B\&B solver on BoxQP instances where $N=\emptyset$. For each instance, we report the number of nodes and the computational time in seconds. When the time limit of 7200 seconds is reached we report the final percentage gap in the column ``Nodes''.}}
\footnotesize
\begin{tabular}{lrr|rr|rr|rr}
\toprule%
& \multicolumn{2}{@{}c@{}}{BARON} & \multicolumn{2}{@{}c@{}}{Gurobi} & \multicolumn{2}{@{}c@{}}{{QuadprogBB}} & \multicolumn{2}{@{}c@{}}{Our B\&B}\\
\cmidrule{2-3}\cmidrule{4-5}\cmidrule{6-7}\cmidrule{8-9}%
Instance	&		Nodes		&	Time	&		Nodes		&	Time	&		Nodes		&	Time	&	Nodes	&	Time	\\
\midrule
s070-025-4	&		18		&	27	&		375		&	2	&		7		&	27	&	1	&	15	\\
s070-025-5	&		16		&	51	&		65		&	1	&		9		&	34	&	1	&	9	\\
s070-025-6	&		24		&	203	&		61		&	2	&		5		&	22	&	1	&	9	\\
s070-050-4	&	(6.92\%)	&	7200	&		4499		&	3658	&		109		&	490	&	1	&	17	\\
s070-050-5	&	(3.92\%)	&	7200	&		83		&	27	&		3		&	18	&	1	&	14	\\
s070-050-6	&	(3.53\%)	&	7200	&		446		&	105	&		21		&	79	&	1	&	16	\\
s070-075-4	&	(10.76\%)	&	7200	&		273		&	369	&		13		&	56	&	1	&	13	\\
s070-075-5	&	(23.84\%)	&	7200	&		610		&	1145	&		431		&	1184	&	1	&	18	\\
s070-075-6	&	(18.89\%)	&	7200	&		359		&	546	&		35		&	110	&	1	&	15	\\
\midrule
s080-025-4	&		2		&	10	&		17		&	1	&		5		&	30	&	1	&	11	\\
s080-025-5	&	(5.07\%)	&	7200	&		6014		&	402	&		47		&	255	&	1	&	21	\\
s080-025-6	&		30		&	761	&		29		&	3	&		13		&	62	&	1	&	17	\\
s080-050-4	&	(4.61\%)	&	7200	&		1339		&	1891	&		47		&	429	&	1	&	20	\\
s080-050-5	&	(4.92\%)	&	7200	&		268		&	315	&		17		&	100	&	1	&	17	\\
s080-050-6	&	(5.16\%)	&	7200	&		358		&	471	&		41		&	248	&	1	&	22	\\
s080-075-4	&	(62.71\%)	&	7200	&		782		&	5786	&		167		&	822	&	4	&	69	\\
s080-075-5	&	(50.0\%)	&	7200	&		303		&	1254	&		437		&	2097	&	1	&	22	\\
s080-075-6	&	(62.13\%)	&	7200	&		392		&	1414	&		33		&	202	&	4	&	89	\\
\midrule
s090-025-4	&		25		&	181	&		19		&	2	&		11		&	93	&	1	&	13	\\
s090-025-5	&	(5.68\%)	&	7200	&	(1.96\%)	&	7200	&		81		&	540	&	1	&	25	\\
s090-025-6	&	(7.49\%)	&	7200	&	(3.41\%)	&	7200	&		209		&	1456	&	1	&	29	\\
s090-050-4	&	(10.97\%)	&	7200	&	(3.00\%)	&	7200	&		135		&	851	&	1	&	27	\\
s090-050-5	&	(14.15\%)	&	7200	&	(0.93\%)	&	7200	&		149		&	973	&	1	&	41	\\
s090-050-6	&	(6.32\%)	&	7200	&		267		&	668	&		17		&	170	&	1	&	17	\\
s090-075-4	&	(103.4\%)	&	7200	&	(2.26\%)	&	7200	&		333		&	2654	&	1	&	58	\\
s090-075-5	&	(148.23\%)	&	7200	&	(8.84\%)	&	7200	&		583		&	4787	&	4	&	70	\\
s090-075-6	&	(97.31\%)	&	7200	&	(1.01\%)	&	7200	&		207		&	1704	&	1	&	29	\\
\midrule
s100-025-4	&	(4.04\%)	&	7200	&		1202		&	485	&		153		&	1450	&	1	&	35	\\
s100-025-5	&	(13.64\%)	&	7200	&	(8.64\%)	&	7200	&	(0.59\%)	&	7200	&	28	&	448	\\
s100-025-6	&	(5.04\%)	&	7200	&	(1.12\%)	&	7200	&		55		&	514	&	1	&	27	\\
s100-050-4	&	(34.12\%)	&	7200	&	(11.17\%)	&	7200	&	(15.67\%)	&	7200	&	55	&	568	\\
s100-050-5	&	(21.97\%)	&	7200	&	(1.38\%)	&	7200	&		117		&	1080	&	1	&	30	\\
s100-050-6	&	(12.82\%)	&	7200	&		179		&	575	&		179		&	1690	&	1	&	32	\\
s100-075-4	&	(132.03\%)	&	7200	&	(4.37\%)	&	7200	&	(4.04\%)	&	7200	&	19	&	541	\\
s100-075-5	&	(141.94\%)	&	7200	&	(4.07\%)	&	7200	&	(0.29\%)	&	7200	&	1	&	59	\\
s100-075-6	&	(129.22\%)	&	7200	&	(7.63\%)	&	7200	&	(10.12\%)	&	7200	&	13	&	238	\\
\midrule
s125-025-4	&	(8.09\%)	&	7200	&	(4.70\%)	&	7200	&		21		&	579	&	1	&	46	\\
s125-025-5	&	(6.03\%)	&	7200	&	(11.75\%)	&	7200	&		49		&	1044	&	1	&	70	\\
s125-025-6	&	(9.72\%)	&	7200	&	(5.75\%)	&	7200	&		433		&	6979	&	1	&	96	\\
s125-050-4	&	(88.03\%)	&	7200	&	(7.08\%)	&	7200	&	(11.92\%)	&	7200	&	70	&	1760	\\
s125-050-5	&	(104.03\%)	&	7200	&	(8.99\%)	&	7200	&	(13.15\%)	&	7200	&	13	&	400	\\
s125-050-6	&	(70.84\%)	&	7200	&	(0.67\%)	&	7200	&		235		&	3878	&	1	&	71	\\
s125-075-4	&	(198.17\%)	&	7200	&	(9.68\%)	&	7200	&	(6.69\%)	&	7200	&	7	&	324	\\
s125-075-5	&	(197.9\%)	&	7200	&	(35.00\%)	&	7200	&	(14.71\%)	&	7200	&	73	&	1924	\\
s125-075-6	&	(192.4\%)	&	7200	&	(18.33\%)	&	7200	&	(7.39\%)	&	7200	&	7	&	288	\\\bottomrule
\end{tabular}
    \label{tab:ternary_results}
\end{table}

\subsection{Experiments over Binary QP problems}\label{sec:binqp}
As a final set of experiments, we tested our approaches over the Binary QP instances and compared them with {{\tt Gurobi} and} the state-of-the-art solver {\tt BiqBin} (see \cite{gusmeroli2022biqbin}). The latter solver strongly relies on a sophisticated parallel implementation based on the Message Passing Interface (MPI) library and it has been optimized to run on High Performance Computing (HPC) systems. While our B\&B search can be parallelized, a more effective parallel implementation of the overall algorithm would require efforts beyond the scope of this paper. Therefore, {to conduct a comparison considering only the algorithmic tools employed}, we evaluate the serial version of {\tt BiqBin} against the serial version of our B\&B with and without multiple variable fixing.
{Note that Binary QP problems can be converted into BoxQP problems as discussed in Section \ref{sec:boxqp}. However, such conversion is not strictly necessary to apply our approach. Indeed, since binary requirements for the variables are already imposed in the problem definition, we could directly apply our approach with $N=\{1,\ldots,n\}$ even if some or all the diagonal elements of matrix {\bf Q} are positive. However, in practice, we observed that there is an advantage in terms of computing times if the Binary QP problem is first converted into a Box QP problem as shown in (\ref{prob:BoxQPLambda}). 
Before presenting the results, we also briefly comment the choice
of vector $\boldsymbol{\lambda}$ in (\ref{prob:BoxQPLambda}).} We need to select vector $\boldsymbol{\lambda}$ in such a way that
$\lambda_i\leq Q_{ii}$. We explored two different definitions for the vector, namely $\boldsymbol{\lambda}=\textrm{diag}({\bf Q})$ and $\boldsymbol{\lambda}=\max_{i} Q_{ii} \mathbf{1}^n$.
The former choice would turn out to be better if local searches were performed via a continuous local solver such as {\tt SNOPT}, because the latter choice induces larger barriers between local minimizers (but for Binary QPs we decided to adopt a discretized local search procedure). On the other hand, the latter choice turned out to be slightly better when included in the B\&B approach. Thus, we adopted it in our experiments. 
{A possible explanation of this fact is that larger values of $\boldsymbol{\lambda}$ do not modify the objective function at the vertices of the box, but increase it in all other parts of the box.}
Concerning the randomized rounding procedure described in Section \ref{sec:standard}, we draw 1000 random samples and we keep the local {minimizer} with the smallest objective value.

The results are displayed in Table \ref{tab:binary_instances_test}. Each instance in the test suite is solved to a final relative optimality gap tolerance of 0.01\% and with a time limit of 7200 seconds. {We also tested the {\tt BARON} solver on these instances. However, we do not report its results since it always reaches the time limit. As already commented in other parts, this is not surprising, since {\tt BARON} is a general-purpose solver.
{\tt Gurobi} is able to solve only small and sparse instances ({\tt be120.3}) within the time limit. Our B\&B solver is competitive with {\tt BiqBin} and
the results confirm that multiple fixing is highly beneficial over the most challenging instances, where a large number of nodes is generated}. The computing times of the version with multiple fixing are often significantly better than those of {\tt BiqBin}. 
{Of course, we stress once again that the performance of {\tt BiqBin} improves considerably if parallelism is allowed. Therefore, we are not claiming that we perform better than {\tt BiqBin}, but just that the algorithmic tools employed in our approach appear to be effective. These results also seem to suggest that a more sophisticated parallel version of our approaches may also lead to competitive results and may be an interesting topic for future research.}   

\begin{figure}[!htbp]
    \centering
    \includegraphics[scale=0.7]{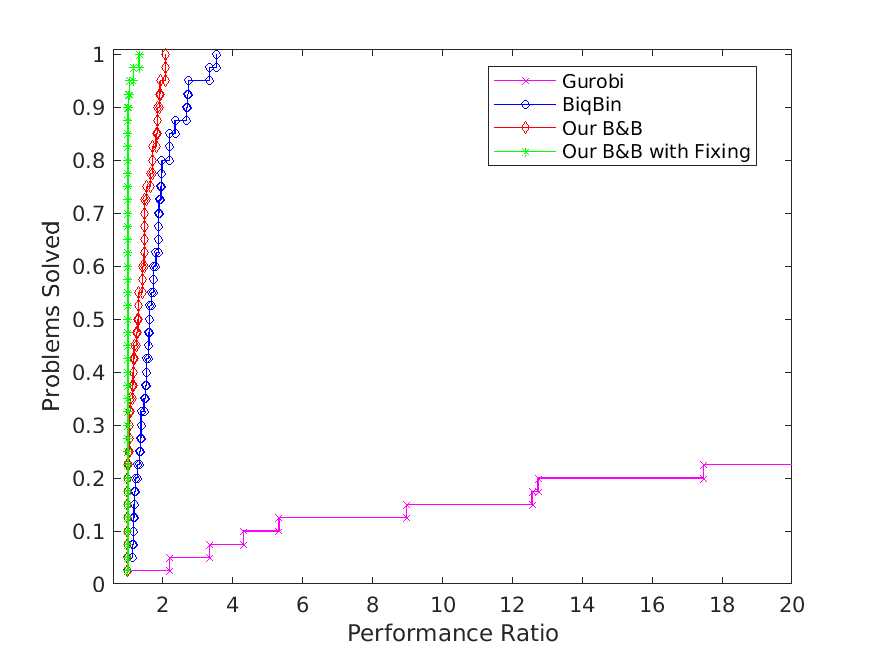}
    \caption{{Computational time performance profile comparing solvers on Binary QP instances ({\tt be$\{n\}$-$\{d\}$-$\{s\}$}, for $n \in \{120, 150\}$, $d \in \{30, 80\}$, and $s \in \{1, 2, 3, 4, 5, 6, 7, 8, 9,10\}$).}}
    \label{fig:perf_binQP}
\end{figure}

{For the tested instances we do not report full statistics on fixing, since they are similar to the ones reported in Table \ref{tab:fix_stat} for the Box QP instances. Note that, when fixing is applied, the reduction of computing times is still significantly larger than the reduction of the number of nodes, confirming that the reduction of time is mostly due to the smaller size of the SDPs solved at nodes of the B\&B tree. }
{In Figure \ref{fig:perf_binQP} the performance profile of the computing time is reported, confirming the effectiveness of our approach.}

\begin{table}[!htbp]
    \centering
    \footnotesize
    \caption{{Comparison between Gurobi, BiqBin our B\&B solver without and with multiple variable fixing on Binary QP instances. Columns ``Nodes'' and ``Time'' indicate the number of nodes and computational time in seconds, respectively; column ``Time Red.'' denotes the percentage reduction in computational time for B\&B with fixing compared to B\&B without fixing. When the time limit of 7200 seconds is reached we report the final percentage gap in the column ``Nodes''.}}
\begin{tabular}{lrr|rr|rr|rrr}
\toprule%
& \multicolumn{2}{@{}c@{}}{Gurobi} &  \multicolumn{2}{@{}c@{}}{BiqBin} & \multicolumn{2}{@{}c@{}}{Our B\&B} & \multicolumn{3}{@{}c@{}}{Our B\&B with Fixing}\\
\cmidrule{2-3}\cmidrule{4-5}\cmidrule{6-7}\cmidrule{8-10}
Instance	&		Nodes		&	Time	&	Nodes	&	Time	&	Nodes	&	Time	&	Nodes	&	Time	&	Time Red. \\	
\midrule
be120.3.1	&		3456		&	403	&	7	&	249	&	3	&	211	&	3	&	184	&	12.80	\\
be120.3.2	&		1264		&	133	&	1	&	49	&	1	&	26	&	1	&	25	&	3.85	\\
be120.3.3	&		3647		&	524	&	1	&	54	&	1	&	30	&	1	&	30	&	0	\\
be120.3.4	&		1738		&	214	&	1	&	74	&	1	&	65	&	1	&	64	&	1.54	\\
be120.3.5	&		3031		&	433	&	1	&	55	&	1	&	34	&	1	&	34	&	0	\\
be120.3.6	&		1662		&	233	&	1	&	38	&	1	&	26	&	1	&	26	&	0	\\
be120.3.7	&		608		&	18	&	1	&	34	&	1	&	21	&	1	&	21	&	0	\\
be120.3.8	&		1158		&	99	&	1	&	27	&	1	&	23	&	1	&	23	&	0	\\
be120.3.9	&		23815		&	6661	&	1	&	82	&	3	&	107	&	3	&	108	&	-0.93	\\
be120.3.10	&		4897		&	716	&	1	&	57	&	1	&	60	&	1	&	60	&	0	\\
\midrule
be120.8.1	&	(12.89\%)	&	7200	&	25	&	752	&	21	&	589	&	21	&	502	&	14.77	\\
be120.8.2	&	(11.02\%)	&	7200	&	13	&	333	&	11	&	343	&	9	&	208	&	39.36	\\
be120.8.3	&	(13.36\%)	&	7200	&	11	&	270	&	7	&	227	&	5	&	155	&	31.72	\\
be120.8.4	&	(11.83\%)	&	7200	&	3	&	139	&	3	&	131	&	3	&	123	&	6.11	\\
be120.8.5	&	(10.49\%)	&	7200	&	1	&	87	&	1	&	73	&	1	&	73	&	0	\\
be120.8.6	&	(14.51\%)	&	7200	&	9	&	293	&	7	&	200	&	3	&	108	&	46	\\
be120.8.7	&	(13.23\%)	&	7200	&	25	&	739	&	15	&	656	&	15	&	315	&	51.98	\\
be120.8.8	&	(16.19\%)	&	7200	&	51	&	1072	&	35	&	1106	&	35	&	779	&	29.57	\\
be120.8.9	&	(15.27\%)	&	7200	&	13	&	428	&	7	&	233	&	7	&	157	&	32.62	\\
be120.8.10	&	(12.20\%)	&	7200	&	15	&	367	&	3	&	263	&	3	&	137	&	47.91	\\
\midrule
be150.3.1	&	(11.92\%)	&	7200	&	11	&	328	&	1	&	93	&	1	&	93	&	0	\\
be150.3.2	&	(13.51\%)	&	7200	&	13	&	337	&	3	&	226	&	3	&	201	&	11.06	\\
be150.3.3	&	(11.17\%)	&	7200	&	1	&	71	&	1	&	38	&	1	&	38	&	0	\\
be150.3.4	&	(10.02\%)	&	7200	&	3	&	264	&	1	&	79	&	1	&	80	&	-1.27	\\
be150.3.5	&	(13.78\%)	&	7200	&	17	&	601	&	7	&	471	&	7	&	322	&	31.63	\\
be150.3.6	&	(16.21\%)	&	7200	&	39	&	1370	&	25	&	1315	&	21	&	895	&	31.94	\\
be150.3.7	&	(13.36\%)	&	7200	&	11	&	385	&	3	&	209	&	3	&	198	&	5.26	\\
be150.3.8	&	(13.98\%)	&	7200	&	63	&	2196	&	33	&	1896	&	33	&	1113	&	41.30	\\
be150.3.9	&	(19.64\%)	&	7200	&	101	&	2382	&	65	&	2071	&	65	&	2072	&	-0.05	\\
be150.3.10	&	(12.88\%)	&	7200	&	65	&	2214	&	37	&	2328	&	33	&	1362	&	41.49	\\
\midrule
be150.8.1	&	(14.02\%)	&	7200	&	49	&	1430	&	29	&	1388	&	27	&	1079	&	22.26	\\
be150.8.2	&	(14.60\%)	&	7200	&	109	&	3937	&	99	&	4643	&	93	&	3253	&	29.94	\\
be150.8.3	&	(14.03\%)	&	7200	&	67	&	2639	&	35	&	2176	&	31	&	1205	&	44.62	\\
be150.8.4	&	(14.62\%)	&	7200	&	21	&	656	&	15	&	648	&	15	&	519	&	19.91	\\
be150.8.5	&	(13.85\%)	&	7200	&	43	&	1482	&	15	&	991	&	13	&	676	&	31.79	\\
be150.8.6	&	(15.12\%)	&	7200	&	41	&	1751	&	31	&	1546	&	31	&	1013	&	34.48	\\
be150.8.7	&	(13.39\%)	&	7200	&	115	&	3345	&	63	&	4043	&	55	&	2185	&	45.96	\\
be150.8.8	&	(14.25\%)	&	7200	&	51	&	1692	&	31	&	1702	&	21	&	898	&	47.24	\\
be150.8.9	&	(16.46\%)	&	7200	&	83	&	3146	&	79	&	4065	&	57	&	1967	&	51.61	\\
be150.8.10	&	(13.99\%)	&	7200	&	43	&	1327	&	23	&	1253	&	21	&	965	&	22.98	\\
\bottomrule
\end{tabular}
    \label{tab:binary_instances_test}
\end{table}

\section{Conclusions}\label{sec:conclusions}
In this paper, we propose an SDP-based branch-and-bound algorithm for solving nonconvex BoxQP problems to global optimality. We first observe that lower bounds from large-scale SDPs with RLT inequalities can be efficiently computed by using first-order SDP solvers. We next observe that adding triangle inequalities in a cutting-plane generation scheme is computationally effective at reducing the optimality gap of BoxQP. We propose a novel fixing strategy allowing to fix multiple variables to 0 or 1 when visiting the B\&B tree. It can be seen as a special case of bound tightening, where the range of variables are reduced to singletons. The most important effect of fixing variables is the size reduction, allowing to solve SDPs over lower dimensional positive semidefinite cones. Computational results in this paper demonstrate the advantage of multiple fixing within a branch-and-bound framework for solving large-scale non-convex BoxQPs. {Possible future developments include:
\begin{itemize} 
\item The exploration of further strategies for multiple variable fixing, aimed both at reducing its cost and at increasing its effectiveness, not only by increasing the number of fixed variables but, e.g., also by selecting variables which, once fixed, allow for a larger increase of the computed lower bounds. A similar idea was exploited in \cite{Sahinidis21} to define the spectral branching;
\item Applying the multiple fixing procedure to more general QP problems, like, e.g., linearly constrained or quadratically constrained (QCQP) problems when binary restrictions are imposed over some or all the variables; for such problems it is not possible to impose binary values of continuous variables $x_i$ such that $Q_{ii}\leq 0$, but the proposed fixing procedure can be applied on binary variables;
\item Development of a solver that tackles larger instances. A potential bottleneck at larger dimensions could be the resolution of SDPs involving all RLTs. As a result, we believe that addressing larger instances may require incorporating both RLTs and triangle inequalities in a cutting-plane procedure;
\item Development of a more advanced parallel version of our approach.
\end{itemize}}

\bibliographystyle{spmpsci}      
\bibliography{abbr, sn-bibliography}

\begin{thebibliography}{10}
\providecommand{\url}[1]{{#1}}
\providecommand{\urlprefix}{URL }
\expandafter\ifx\csname urlstyle\endcsname\relax
  \providecommand{\doi}[1]{DOI~\discretionary{}{}{}#1}\else
  \providecommand{\doi}{DOI~\discretionary{}{}{}\begingroup \urlstyle{rm}\Url}\fi

\bibitem{alizadeh1995interior}
Alizadeh, F.: Interior point methods in semidefinite programming with applications to combinatorial optimization.
\newblock SIAM J. Optim. \textbf{5}(1), 13--51 (1995)

\bibitem{anstreicher2009semidefinite}
Anstreicher, K.M.: Semidefinite programming versus the reformulation-linearization technique for nonconvex quadratically constrained quadratic programming.
\newblock J. Global Optim. \textbf{43}(2), 471--484 (2009)

\bibitem{Bomze02}
Bomze, I.M., De~Klerk, E.: Solving standard quadratic optimization problems via linear, semidefinite and copositive programming.
\newblock J. Global Optim. \textbf{24}(2), 163--185 (2002)

\bibitem{Bonami18}
Bonami, P., G{\"u}nl{\"u}k, O., Linderoth, J.: Globally solving nonconvex quadratic programming problems with box constraints via integer programming methods.
\newblock Math. Program. Comput. \textbf{10}(3), 333--382 (2018)

\bibitem{Burer09}
Burer, S.: On the copositive representation of binary and continuous nonconvex quadratic programs.
\newblock Math. Program. \textbf{120}(2), 479--495 (2009)

\bibitem{Burer10}
Burer, S.: Optimizing a polyhedral-semidefinite relaxation of completely positive programs.
\newblock Math. Program. Comput. \textbf{2}(1), 1--19 (2010)

\bibitem{burer2009nonconvex}
Burer, S., Letchford, A.N.: On nonconvex quadratic programming with box constraints.
\newblock SIAM J. Optim. \textbf{20}(2), 1073--1089 (2009)

\bibitem{Burer09a}
Burer, S., Vandenbussche, D.: Globally solving box-constrained nonconvex quadratic programs with semidefinite-based finite branch-and-bound.
\newblock Comput. Optim. Appl. \textbf{43}(2), 181--195 (2009)

\bibitem{Caprara10}
Caprara, A., Locatelli, M.: Global optimization problems and domain reduction strategies.
\newblock Math. Program. \textbf{125}(1), 123--137 (2010)

\bibitem{cerulli2021improving}
Cerulli, M., De~Santis, M., Gaar, E., Wiegele, A.: Improving admms for solving doubly nonnegative programs through dual factorization.
\newblock 4OR \textbf{19}(3), 415--448 (2021)

\bibitem{Chen12}
Chen, J., Burer, S.: Globally solving nonconvex quadratic programming problems via completely positive programming.
\newblock Math. Program. Comput. \textbf{4}(1), 33--52 (2012)

\bibitem{dolan2002benchmarking}
Dolan, E.D., Mor{\'e}, J.J.: Benchmarking optimization software with performance profiles.
\newblock Math. Program. \textbf{91}(2), 201--213 (2002)

\bibitem{elloumi2019global}
Elloumi, S., Lambert, A.: Global solution of non-convex quadratically constrained quadratic programs.
\newblock Optim. Methods Softw. \textbf{34}(1), 98--114 (2019)

\bibitem{galli2018binarisation}
Galli, L., Letchford, A.N.: A binarisation heuristic for non-convex quadratic programming with box constraints.
\newblock Oper. Res. Lett. \textbf{46}(5), 529--533 (2018)

\bibitem{gill2005snopt}
Gill, P.E., Murray, W., Saunders, M.A.: Snopt: An sqp algorithm for large-scale constrained optimization.
\newblock SIAM Review \textbf{47}(1), 99--131 (2005)

\bibitem{Gleixner17}
Gleixner, A.M., Berthold, T., M{\"u}ller, B., Weltge, S.: Three enhancements for optimization-based bound tightening.
\newblock J. Global Optim. \textbf{67}(4), 731--757 (2017)

\bibitem{goemans1995improved}
Goemans, M.X., Williamson, D.P.: Improved approximation algorithms for maximum cut and satisfiability problems using semidefinite programming.
\newblock J. ACM \textbf{42}(6), 1115--1145 (1995)

\bibitem{Gondzio18}
Gondzio, J., Y{\i}ld{\i}r{\i}m, E.A.: Global solutions of nonconvex standard quadratic programs via mixed integer linear programming reformulations.
\newblock J. Global Optim. \textbf{81}(2), 293--321 (2021)

\bibitem{gusmeroli2022biqbin}
Gusmeroli, N., Hrga, T., Lu{\v{z}}ar, B., Povh, J., Siebenhofer, M., Wiegele, A.: Biqbin: a parallel branch-and-bound solver for binary quadratic problems with linear constraints.
\newblock ACM Trans. Math. Softw. \textbf{48}(2), 1--31 (2022)

\bibitem{Hansen93}
Hansen, P., Jaumard, B., Ruiz, M., Xiong, J.: Global minimization of indefinite quadratic functions subject to box constraints.
\newblock Nav. Res. Logist. \textbf{40}(3), 373--392 (1993)

\bibitem{jansson2008rigorous}
Jansson, C., Chaykin, D., Keil, C.: Rigorous error bounds for the optimal value in semidefinite programming.
\newblock SIAM J. Numer. Anal. \textbf{46}(1), 180--200 (2008)

\bibitem{khajavirad2018hybrid}
Khajavirad, A., Sahinidis, N.V.: A hybrid lp/nlp paradigm for global optimization relaxations.
\newblock Math. Program. Comput. \textbf{10}(3), 383--421 (2018)

\bibitem{Liuzzi19}
Liuzzi, G., Locatelli, M., Piccialli, V.: A new branch-and-bound algorithm for standard quadratic programming problems.
\newblock Optim. Methods Softw. \textbf{34}(1), 79--97 (2019)

\bibitem{Liuzzi22}
Liuzzi, G., Locatelli, M., Piccialli, V.: A computational study on qp problems with general linear constraints.
\newblock Optim. Lett. \textbf{16}(6), 1633--1647 (2022)

\bibitem{Liuzzi20}
Liuzzi, G., Locatelli, M., Piccialli, V., Rass, S.: Computing mixed strategies equilibria in presence of switching costs by the solution of nonconvex {QP} problems.
\newblock Comput. Optim. Appl. \textbf{79}(3), 561--599 (2021)

\bibitem{mccormick1976computability}
McCormick, G.P.: Computability of global solutions to factorable nonconvex programs: Part i—convex underestimating problems.
\newblock Math. Program. \textbf{10}(1), 147--175 (1976)

\bibitem{Motzkin65}
Motzkin, T.S., Straus, E.G.: Maxima for graphs and a new proof of a theorem of tur{\'a}n.
\newblock Can. J. Math. \textbf{17}, 533--540 (1965)

\bibitem{Sahinidis21}
Nohra, C.J., Raghunathan, A.U., Sahinidis, N.: Spectral relaxations and branching strategies for global optimization of mixed-integer quadratic programs.
\newblock SIAM J. Optim. \textbf{31}(1), 142--171 (2021)

\bibitem{padberg1989boolean}
Padberg, M.: The boolean quadric polytope: some characteristics, facets and relatives.
\newblock Math. Program. \textbf{45}(1), 139--172 (1989)

\bibitem{park2018semidefinite}
Park, J., Boyd, S.: A semidefinite programming method for integer convex quadratic minimization.
\newblock Optim. Lett. \textbf{12}(3), 499--518 (2018)

\bibitem{piccialli2022b}
Piccialli, V., {Russo Russo}, A., Sudoso, A.M.: An exact algorithm for semi-supervised minimum sum-of-squares clustering.
\newblock Comput. Oper. Res. \textbf{147}, 105958 (2022)

\bibitem{piccialli2022c}
Piccialli, V., Sudoso, A.M.: Global optimization for cardinality-constrained minimum sum-of-squares clustering via semidefinite programming.
\newblock Math. Program. pp. 1--35 (2023)

\bibitem{piccialli2022a}
Piccialli, V., Sudoso, A.M., Wiegele, A.: {SOS-SDP}: An exact solver for minimum sum-of-squares clustering.
\newblock INFORMS J. Comput. \textbf{34}(4), 2144--2162 (2022)

\bibitem{Sahinidis96}
Sahinidis, N.V.: {BARON}: A general purpose global optimization software package.
\newblock J. Global Optim. \textbf{8}(2), 201--205 (1996)

\bibitem{sahinidis1996baron}
Sahinidis, N.V.: Baron: A general purpose global optimization software package.
\newblock J. Global Optim. \textbf{8}, 201--205 (1996)

\bibitem{sherali1995reformulation}
Sherali, H.D., Tuncbilek, C.H.: A reformulation-convexification approach for solving nonconvex quadratic programming problems.
\newblock J. Global Optim. \textbf{7}(1), 1--31 (1995)

\bibitem{shor1987quadratic}
Shor, N.Z.: Quadratic optimization problems.
\newblock Sov. J. Comput. Syst. Sci. \textbf{25}, 1--11 (1987)

\bibitem{sun2015convergent}
Sun, D., Toh, K.C., Yang, L.: A convergent 3-block semiproximal alternating direction method of multipliers for conic programming with 4-type constraints.
\newblock SIAM J. Optim. \textbf{25}(2), 882--915 (2015)

\bibitem{sun2020sdpnal+}
Sun, D., Toh, K.C., Yuan, Y., Zhao, X.Y.: {SDPNAL+}: A {M}atlab software for semidefinite programming with bound constraints (version 1.0).
\newblock Optim. Methods Softw. \textbf{35}(1), 87--115 (2020)

\bibitem{Tawarmalani04}
Tawarmalani, M., Sahinidis, N.V.: Global optimization of mixed-integer nonlinear programs: A theoretical and computational study.
\newblock Math. Program. \textbf{99}(3), 563--591 (2004)

\bibitem{wiegele2007biq}
Wiegele, A.: Biq mac library (2007).
\newblock \urlprefix\url{https://biqmac.aau.at/biqmaclib.html}

\bibitem{Xia20}
Xia, W., Vera, J.C., Zuluaga, L.F.: Globally solving nonconvex quadratic programs via linear integer programming techniques.
\newblock INFORMS J. Comput. \textbf{32}(1), 40--56 (2020)

\bibitem{Yajima1998}
Yajima, Y., Fujie, T.: A polyhedral approach for nonconvex quadratic programming problems with box constraints.
\newblock J. Global Optim. \textbf{13}(2), 151--170 (1998)

\bibitem{yang2015sdpnal}
Yang, L., Sun, D., Toh, K.C.: {SDPNAL+}: a majorized semismooth {N}ewton-{CG} augmented {L}agrangian method for semidefinite programming with nonnegative constraints.
\newblock Math. Program. Comput. \textbf{7}(3), 331--366 (2015)

\end{thebibliography}


\end{document}